\documentclass[11pt]{article}

\usepackage{amsmath, amssymb, amscd}

\def\eqref#1{(\ref{#1})}
\newcommand{\goth}{\frak}

\newcommand{\arrow}{{\:\longrightarrow\:}}

\newcommand{\C}{{\Bbb C}}
\newcommand{\R}{{\Bbb R}}

\renewcommand{\H}{{\Bbb H}}
\newcommand{\6}{\partial}
\newcommand{\1}{\sqrt{-1}\:}
\newcommand{\restrict}[1]{{\left|_{{\phantom{|}\!\!}_{#1}}\right.}}

\renewcommand{\c}[1]{{\cal #1}}
\newcommand{\calo}{{\cal O}}

\renewcommand{\tilde}{\widetilde}
\renewcommand{\bar}{\overline}
\renewcommand{\phi}{\varphi}
\renewcommand{\epsilon}{\varepsilon}
\renewcommand{\geq}{\geqslant}
\renewcommand{\leq}{\leqslant}

\newcommand{\End}{\operatorname{End}}

\newcommand{\Id}{\operatorname{Id}}
\newcommand{\id}{\operatorname{\text{\sf id}}}
\newcommand{\Vol}{\operatorname{Vol}}
\newcommand{\Hom}{\operatorname{Hom}}

\newcommand{\codim}{\operatorname{codim}}

\newcommand{\rk}{\operatorname{rk}}

\newcommand{\Tr}{\operatorname{Tr}}

\renewcommand{\Re}{\operatorname{\sf Re}}

\newcommand{\comment}[1]{{}}

\def\blacksquare{\hbox{\vrule width 4pt height 4pt depth 0pt}}
\def\endproof{\blacksquare}

\makeatletter

\@ifundefined{Bbb}
     {\newcommand{\Bbb}[1]{{\mathbb #1}}}%
{}

\newcommand{\ps@verbit}{%
  \renewcommand{\@oddhead}{%
          \scriptsize
          {Hyperholomorphic connections on sheaves}
          \hfil\tiny {M. Verbitsky, \ \ \ \ July 16, 2001}}
  \renewcommand{\@evenhead}{\@oddhead}
  \renewcommand{\@oddfoot}{\hfil\thepage\hfil}
  \renewcommand{\@evenfoot}{\@oddfoot}}
 
\newcounter{Mycounter}[section]
\newcounter{lemma}[section]
\renewcommand{\thelemma}{{Lemma \thesection.\arabic{lemma}}}
\newcommand{\lemma}{%
     \setcounter{lemma}{\value{Mycounter}}
     \refstepcounter{lemma}
     \stepcounter{Mycounter}
     {\bf \thelemma:\ }}

\newcounter{claim}[section]
\renewcommand{\theclaim}{{Claim \thesection.\arabic{claim}}}
\newcommand{\claim}{%
     \setcounter{claim}{\value{Mycounter}}
     \refstepcounter{claim}
     \stepcounter{Mycounter}
     {\bf \theclaim:\ }}

\newcounter{sublemma}[section]
\newcounter{corollary}[section]
\renewcommand{\thecorollary}{{Corollary \thesection.\arabic{corollary}}}
\newcommand{\corollary}{%
     \setcounter{corollary}{\value{Mycounter}}
     \refstepcounter{corollary}
     \stepcounter{Mycounter}
     {\bf \thecorollary:\ }}

\newcounter{theorem}[section]
\renewcommand{\thetheorem}{{Theorem \thesection.\arabic{theorem}}}
\newcommand{\theorem}{%
     \setcounter{theorem}{\value{Mycounter}}
     \refstepcounter{theorem}
     \stepcounter{Mycounter}
     {\bf \thetheorem:\ }}

\newcounter{conjecture}[section]
\renewcommand{\theconjecture}{{Conjecture \thesection.\arabic{conjecture}}}
\newcommand{\conjecture}{%
     \setcounter{conjecture}{\value{Mycounter}}
     \refstepcounter{conjecture}
     \stepcounter{Mycounter}
     {\bf \theconjecture:\ }}

\newcounter{proposition}[section]
\renewcommand{\theproposition}
       {{Proposition \thesection.\arabic{proposition}}}
\newcommand{\proposition}{%
     \setcounter{proposition}{\value{Mycounter}}
     \refstepcounter{proposition}
     \stepcounter{Mycounter}
     {\bf \theproposition:\ }}

\newcounter{definition}[section]
\renewcommand{\thedefinition}
       {{Definition~\thesection.\arabic{definition}}}
\newcommand{\definition}{%
     \setcounter{definition}{\value{Mycounter}}
     \refstepcounter{definition}
     \stepcounter{Mycounter}
     {\bf \thedefinition:\ }}

\newcounter{example}[section]
\renewcommand{\theexample}{{Example \thesection.\arabic{example}}}
\newcommand{\example}{%
     \setcounter{example}{\value{Mycounter}}
     \refstepcounter{example}
     \stepcounter{Mycounter}
     {\bf \theexample:\ }}

\newcounter{remark}[section]
\renewcommand{\theremark}{{Remark \thesection.\arabic{remark}}}
\newcommand{\remark}{%
     \setcounter{remark}{\value{Mycounter}}
     \refstepcounter{remark}
     \stepcounter{Mycounter}
     {\bf \theremark:\ }}

\newcounter{problem}[section]
\newcounter{question}[section]
\@addtoreset{equation}{section}
\@addtoreset{footnote}{section}
\makeatother

\begin{document}

\begin{center}
{\LARGE\bf
Hyperholomorphic connections   \\[3mm] on coherent  sheaves 
and stability
}
\\[4mm]
Misha Verbitsky,\footnote{ The author is
supported by the RFBR grant 09-01-00242-a,
 RFBR grant 10-01-93113-NCNIL-a, Science Foundation of 
the SU-HSE award No. 10-09-0015, and AG Laboratory 
SU-HSE, RF government grant, ag. 11.G34.31.0023.}
\\[4mm]
{\tt verbit@thelema.dnttm.ru, verbit@mccme.ru}
\end{center}

{\small 
\hspace{0.15\linewidth}
\begin{minipage}[t]{0.75\linewidth}
{\bf Abstract} \\
Let $M$ be a hyperk\"ahler manifold, 
and $F$ a reflexive sheaf on $M$.
Assume that $F$ (outside of singularities)
admits a connection $\nabla$ with a curvature
$\Theta$ which is invariant under the standard
$SU(2)$-action on 2-forms. If $\Theta$ is square-integrable,
such sheaf is called {\bf hyperholomorphic}.
Hyperholomorphic sheaves were studied at great
length in \cite{_V:Hyperholo_sheaves_}. Such sheaves
are stable and their singular sets are hyperk\"ahler
subvarieties in $M$. In the present paper, we study sheaves 
admitting a connection with $SU(2)$-invariant curvature 
which is not necessary $L^2$-integrable. We show that such
sheaves are polystable.
\end{minipage}
}

{
\small
\tableofcontents
}

\section{Introduction}


Yang-Mills theory of holomorphic vector bundles
is one of the most spectacular successes of modern algebraic
geometry. Developed by Narasimhan-Seshadri, Kobayashi,
Hitchin, Donaldson, Uhlenbeck-Yau and others,
this theory proved to be very fruitful
in the study of stability and the 
modular properties 
of holomorphic vector bundles.
Bogomolov-Miyaoka-Yau inequality and 
Uhlenbeck-Yau theorem were used by Carlos Simpson in his
groundbreaking works on variations of Hodge structures 
and flat bundles (\cite{_Simpson_}). Later, it was shown
(\cite{_Verbitsky:Hyperholo_bundles_}, \cite{_NHYM_})
that Yang-Mills approach is also useful in hyperk\"ahler
geometry and lends itself to an extensive study of stable
bundles, their modular and twistor properties. 

{}From algebraic point of view, a coherent sheaf is much more
natural kind of object than a holomorphic vector bundle.
This precipitates the extreme importance 
of Bando-Siu theory \cite{_Bando_Siu_} 
which extends Yang-Mills geometry
to coherent sheaves.

In the present paper, we study the ramifications of 
Bando-Siu theory, for hyperk\"ahler manifolds.

\subsection{Yang-Mills geometry and Bando-Siu theorem}

Let $M$ be a compact K\"ahler manifold, and $\omega$ its K\"ahler
form. Consider the standard Hodge operator $L$ on differential forms
which multiplies a form by $\omega$. Let $\Lambda$ be the Hermitian 
adjoint operator. 

Let $B$ be a Hermitian holomorphic vector bundle,
and \[ \Theta\in \Lambda^{1,1}(M, \End(B))\] its curvature,
considered as a (1,1)-form on $M$ with coefficients in
$\End(B)$. By definition, $\Lambda \Theta$ is a 
smooth section of $\End(B)$.
The bundle $B$ is called {\bf Yang-Mills},
or {\bf Hermitian-Einstein}, if 
$\Lambda\Theta$ is constant times the unit section of
$\End(B)$. 

For a definition of stability, see 
Subsection \ref{_sta_bu_and_YM_Subsection_}. Throughout
this paper, stability is understood in the sense of
Mumford-Takemoto.

Let $B$ be a holomorphic bundle which cannot be decomposed
onto a direct sum of non-trivial holomorphic bundles
(such bundles are called {\bf indecomposable}). 
By Uhlenbeck-Yau theorem,
$B$ admits a
Yang-Mills connection if and only if $B$ is stable;
if exists, such a connection is unique 
(\cite{_Uhle_Yau_}). This result allows one to deal
with the moduli of stable holomorphic vector bundles
in efficient and straightforward manner. 

S. Bando and Y.-T. Siu (\cite{_Bando_Siu_})
extended the results of Uhlenbeck-Yau to coherent sheaves.

\hfill

\definition\label{_admissi_Intro_Definition_}
\cite{_Bando_Siu_} 
Let $F$ be a coherent sheaf on $M$ and
$\nabla$ a Hermitian connection
on $F$ defined outside of its singularities.
Denote by $\Theta$ the curvature of $\nabla$.
Then $\nabla$ is called {\bf admissible}
if the following holds
\begin{description}
\item[(i)] $\Lambda \Theta\in \End(F)$ is uniformly bounded
\item[(ii)] $|\Theta|^2$ is integrable on $M$. 
\end{description}

\hfill

Any torsion-free 
coherent sheaf admits an  admissible connection.
An admissible connection can be extended over the
place where $F$ is smooth. Moreover, if a bundle $B$ 
on $M\backslash Z$, $\codim_\C Z\geq 2$ is equipped
with an admissible connection, then $B$ can be extended to
a coherent sheaf on $M$. 

Therefore, the notion of a coherent sheaf can be adequately replaced
by the notion of an admissible Hermitian holomorphic bundle
on $M\backslash Z$, $\codim_\C Z\geq 2$.

A version of Uhlenbeck-Yau theorem exists for 
coherent sheaves (\ref{_UY_for_shea_Theorem_});
given a torsion-free coherent sheaf $F$, $F$ admits an admissible
Yang-Mills connection $\nabla$ if and only if $F$ is polystable.

The following conjecture deals with Yang-Mills connections
which are {\it not} admissible.

\hfill

\conjecture\label{_admiss_YM_Conjecture_}
Let $M$ be a compact K\"ahler manifold, $F$ a torsion-free
coherent sheaf on $M$ with singularities in codimension
at least 3, $\nabla$ a Hermitian connection on $F$
defined outside of its singularities, and $\Theta$ its
curvature. Assume that
$\Lambda(\Theta)=0$. Then $F$ 
can be extended to a stable sheaf on $M$.

\hfill

This conjecture is motivated by the following heuristic argument.

Denote by $\omega$ the K\"ahler form on $M$, and 
let $n:= \dim_\C M$.
By Hodge-Riemann relations,

\begin{equation}\label{_H-R_for_curva_Equation_}
\Tr(\Theta\wedge \Theta) \wedge \omega^{n-2} = c | \Theta|^2 \Vol(M)
\end{equation}
where $c$ is a positive rational 
constant (this equality is true pointwise,
assuming that $\Lambda \Theta =0$). This 
equality is used in \cite{_Simpson_} to deduce
the Bogomolov-Miyaoka-Yau inequality
from the Uhlenbeck-Yau theorem.

By Gauss-Bonnet formula,
the cohomology class of $\Tr(\Theta\wedge \Theta)$ can be expressed
via $c_1(F)$, $c_2(F)$:
\[ 
  \frac{\1}{{2\pi}^2}\Tr(\Theta\wedge \Theta) = 2 c_2(F) - \frac{n-1}{n} c_1(F).
\]
Therefore, the integral
\begin{equation}\label{_integral_c_2_Equation_}
   \int_M\Tr(\Theta\wedge \Theta) \wedge \omega^{n-2}    
\end{equation}
``must have'' cohomological meaning (we write ``must have'' to indicate
here the element of speculation). 

If indeed the integral \eqref{_integral_c_2_Equation_}
is expressed via cohomology, it is finite, and by 
\eqref{_H-R_for_curva_Equation_} 
the curvature $\Theta$ is square-integrable. 

In this paper, we study \ref{_admiss_YM_Conjecture_}
when $M$ is a hyperk\"ahler manifold, and $\nabla$ is
a hyperholomorphic connection (see \ref{_hyperho_conne_Definition_} 
for a definition and further discussion of the notion 
of hyperholomorphic bundle).

\subsection{Hyperk\"ahler and hypercomplex manifolds}

A hypercomplex manifold is a manifold 
equipped with an action of quaternion algebra in its
tangent bundle $TM$, such that for any
quaternion $L, L^2=-1$, the corresponding operator
on $TM$ defines an integrable structure on $M$.
If, in addition, $M$ is Riemannian, and $(M, L)$
is K\"ahler for any quaternion $L, L^2=-1$,
then $M$ is called {\bf hyperk\"ahler}.

A hyperk\"ahler manifold is equipped with a natural action of
the group $SU(2)$ on $TM$. By multiplicativity, we may extend
this action to all tensor powers of $TM$. In particular, $SU(2)$
acts on the space of differential forms on $M$.

This action bears a deep geometric meaning encompassing the
Hodge decomposition on $M$ (see 
\ref{_SU(2)_inva_type_p,p_Lemma_} and its proof).
Moreover, the group $SU(2)$ preserves
the Laplace operator, and henceforth acts 
on the cohomology of $M$
(see e.g. \cite{_Verbitsky:Symplectic_II_}).

Let $\eta$ be an $SU(2)$-invariant 2-form on $M$.
An elementary linear-\-al\-geb\-raic calculation implies
that $\Lambda \eta=0$ (\ref{_Lambda_of_inva_forms_zero_Lemma_}). 
Given a Hermitian vector bundle $B$ with $SU(2)$-invariant
curvature 
\[ \Theta \in \Lambda^2(M)_{SU(2)-inv} \otimes \End(B),
\]
we find that $\Lambda\Theta=0$, and therefore $B$ is Yang-Mills.

Such bundles are called {\bf hyperholomorphic}.
The theory of hyperholomorphic bundles, developed in 
\cite{_Verbitsky:Hyperholo_bundles_}, turns out
to be quite useful in hyperk\"ahler geometry,
by the following reasons.

\begin{equation}\label{_properties_hh_Equation_}
\text{\begin{minipage}[t]{0.7\linewidth}
\begin{description}
\item[(i)] For an arbitrary holomorphic vector bundle,
a hyperholomorphic connection is Yang-Mills, and therefore
unique.
\item[(ii)] An $SU(2)$-invariant form is of the Hodge
type $(1,1)$ with respect to any  complex structure
$L\in \H$, $L^2=-1$ induced by the quaternionic action on $M$
(\ref{_SU(2)_inva_type_p,p_Lemma_}). By Newlander-Nirenberg
integrability theorem
(\ref{_Newle_Nie_for_bu_Theorem_}), a hyperholomorphic bundle
is holomorphic with respect to $I, J, K\in \H$. 
The converse is also true (\ref{_hyperho_conne_Definition_}). 
\end{description}
\end{minipage}}
\end{equation}
\begin{equation*}
\text{\begin{minipage}[t]{0.7\linewidth}
\begin{description}
\item[(iii)] The moduli of hyperholomorphic bundles
are hyperk\"ahler (possibly singular) varieties.
A normalization of such variety is smooth and
hyperk\"ahler (\cite{_Verbitsky:hypercomple_}).
\item[(iv)]  Let $L\in \H$, $L^2=-1$ 
be a complex structure
induced by the quaternionic action.
Consider a stable holomorphic bundle on the K\"ahler
manifold $(M, L)$. Then $B$ admits a hyperholomorphic
connection if and only if the Chern classes
$c_1(B)$, $c_2(B)$ are $SU(2)$-invariant
(\ref{_inva_then_hyperho_Theorem_}).
\end{description}
\end{minipage}}
\end{equation*}
\begin{equation*}
\text{\begin{minipage}[t]{0.7\linewidth}
\begin{description}
 \item[(v)] Moreover, if $L\in \H$, $L^2=-1$  is generic,
and $B$ is a stable holomorphic bundle
on $(M,L)$, then $B$ is hyperholomorphic.
\end{description}
\end{minipage}}
\end{equation*}

Using the results of Bando-Siu, we can extend the notion
of hyperholomorphic connection to coherent sheaves
(\cite{_V:Hyperholo_sheaves_}). 

\hfill

\definition\label{_admissi_hy_Definition_}
Let $M$ be a hyperk\"ahler manifold,
and $F$ a reflexive\footnote{A torsion-free coherent sheaf
is called {\bf reflexive} if the natural monomorphism 
\[ F\arrow F^{**}:=\Hom(\Hom(F, \calo), \calo)\] 
is an isomorphism. A sheaf $F^{**}$ is always reflexive.
The natural functor $F\arrow F^{**}$ is called
{\bf the reflexization}. For more details on reflexive
sheaves, see Subsection \ref{_stable_refle_YM_Subsection_}
and \cite{_OSS_}.}
 coherent sheaf on 
the K\"ahler manifold $(M,I)$. Consider 
an admissible (in the sense of \ref{_admissi_Intro_Definition_})
Hermitian connection $\nabla$ 
on $F$. The $\nabla$ is called 
{\bf admissible hyperholomorphic}, if its curvature
is $SU(2)$-invariant. A stable reflexive sheaf is called
{\bf stable hyperholomorphic} if it admits an 
admissible hyperholomorphic connection.

\hfill

The statements (i)-(ii) and (iv)-(v) 
of \eqref{_properties_hh_Equation_}
hold true for hyperholomorphic sheaves. In addition to this,
a hyperholomorphic sheaf with isolated singularities
can be desingularized with a single blow-up
(\cite{_V:Hyperholo_sheaves_}).

\hfill

In examples, one often obtains coherent sheaves with
Hermitian structure outside of singularities. For instance,
a direct image of a Hermitian vector bundle
is a complex of sheaves with cohomology equipped
with the natural (Weil-Peterson) metrics.
If we work in hyperk\"ahler geometry, the
corresponding Hermitian connection is quite 
often hyperholomorphic outside of singularities
(\cite{_BBR_3}). However, the admissibility condition
is rather tricky. In fact, we were unable to show 
in full generality that a sheaf with a connection
and an $SU(2)$-invariant curvature is admissible.

However, the following assertion is sufficient for most
purposes.

\hfill

\theorem\label{_main_theo_intro_Theorem_}
Let $M$ be a compact hyperk\"ahler manifold, $I$ an induced complex
structure, and $F$ a reflexive sheaf on $(M,I)$
which cannot be decomposed onto a direct sum of non-trivial
coherent sheaves.\footnote{Such sheaves are called indecomposable.} 
Assume that $F$ is equipped with a Hermitian connection $\nabla$ defined 
outside of the singular set of $F$, and the curvature of $\nabla$
is $SU(2)$-invariant. Then $F$ is stable.

\hfill

{\bf Proof:}
This is \ref{_wea_hh_hh_Theorem_}. \endproof

\subsection{Contents}

This paper has the following structure.

\hfill

$\bullet$ The present Introduction is independent from the rest of this
paper. 

\hfill

$\bullet$ In Sections \ref{_hyperka_Section_}-\ref{_hyperho_Section_}
we give preliminary definitions and state basic results
about the geometry of hyperk\"ahler manifolds and stable
bundles. We follow \cite{_Besse:Einst_Manifo_}, 
\cite{_Uhle_Yau_} and \cite{_Verbitsky:Hyperholo_bundles_}.

\hfill

$\bullet$ In Subsections 
\ref{_stable_refle_YM_Subsection_}-\ref{_hyperholo_shea_def_Subsection_},
we give an exposition of the theory of Bando-Siu and its
applications to the hyperk\"ahler geometry. We also 
give a definition of reflexive sheaves and list some
of their properties.

\hfill

$\bullet$ In Subsection \ref{_weak_hype_Subsection_},
we state the main conjecture motivating our research
(\ref{_wea_hh_conn_hh_Conjecture_}). 
It is conjectured that, on any 
hyperk\"ahler manifold, a hyperholomorphic connection on a
reflexive sheaf (defined everywhere 
outside of singularities) has
square-integrable curvature.

We also state our main result (\ref{_wea_hh_hh_Theorem_})
which was explained earlier in this Introduction 
(\ref{_main_theo_intro_Theorem_}). 

\hfill

$\bullet$ In Section \ref{_posi_c_2_Section_} we work with
positive $(p,p)$-forms and their singularities. We state
an important lemma of Sibony, motivating \ref{_wea_hh_conn_hh_Conjecture_}.
Let $\eta$ be a positive closed $(p,p)$-form
with singularities in codimension at least $p+1$.
Then $\eta$ is $L^1$-integrable. This is used to prove 
\ref{_wea_hh_conn_hh_Conjecture_} in case of a sheaf
with isolated singularities.

As an intermediate result, we obtain the following proposition,
which is quite useful in itself
(\ref{_c_2_wedge_omega^n-3_positive_Lemma_}).
Let $B$ be a hyperholomorphic bundle on a hyperk\"ahler manifold
$M$, $\dim_\H M =n >1$,  $\Theta$ its curvature. Denote by
$\omega_I$ the K\"ahler form of $(M,I)$.
Consider the closed 4-form
\[ r_2:= \frac{\1}{{2\pi}^2}\Tr(\Theta\wedge \Theta)
\]
representing (by Gauss-Bonnet) the cohomology class 
$2 c_2(B) - \frac{n-1}{n} c_1(B)^2$. 
Then the $(2n-1,2n-1)$-form
$r_2\wedge \omega_I^{2n-3}$ is positive.

\hfill

$\bullet$ In Section \ref{_SU(2)_inv_Section_}
we study the first Chern class of a reflexive sheaf $F$ admitting a 
hyperholomorphic connection outside of singularities. We  
show that $c_1(F)$ is $SU(2)$-invariant. 

\hfill

$\bullet$ In Section \ref{_positi_hype_Section_}, 
we study the singularities of positive forms on a hyperk\"ahler
manifold. Consider a closed 2-form $\eta$ which is smooth on
$M\backslash Z$, where $codim_\R Z \geq 6$. 
Since $H^2(M\backslash Z) = H^2(M)$, we may consider
the cohomology class $[\eta]$ as an element in $H^2(M)$.
We define the degree $\deg_I$ of $[\eta]$ as follows
\[ \deg_I(\eta) := \int_M [\eta]\wedge \omega_I^{n-1},
\]
where $\omega_I$ is the K\"ahler form of $(M,I)$.

Assume that $\eta$ is a sum of a positive form $\eta_+$ 
and an $SU(2)$-invariant form. We show that $\deg_I[\eta]\geq 0$,
and if $\deg_I[\eta]=0$, then $\eta_+=0$.

\hfill

$\bullet$ In Section \ref{_integra_K-posit_Section_},
we prove $L^1$-integrability 
of a $\6_K$-closed form $\eta^{2,0}_K\in \Lambda^{2,0}_K(M\backslash Z)$, 
where $I(\eta^{2,0}_K) = \bar \eta^{2,0}_K$, assuming that
$\Re\eta(z, \bar z)$ 
is non-negative for all $z\in T^{1,0}(M,I)$. This is essentially
a hyperk\"ahler version of two classical results
from complex analysis - Sibony's lemma and
Skoda-El Mir theorem. This result is used in Section \ref{_positi_hype_Section_}
to show that certain closed forms with singularities
represent cohomology classes 
of positive degree.

In the earlier versions of this paper this result was  
proven by a straightforward
argument based on slicing, in the same way as one proves 
the $L^1$-integrability of a positive closed $(p,p)$ form
with singularities in $\codim >2p$ (\cite{_Sibony_}).
To use slicing, one needs to approximate a hyperk\"ahler 
manifold by a flat one, which leads to complicated estimates. 
Now these difficulties are avoided.
In the latest version (starting from 2008), a coordinate-free
approach to Sibony's lemma was used, based on the
recent advances in the theory of $\omega^q$-plurisubharmonic
functions (\cite{_Verbitsky:omega-psh_}, \cite{_Verbitsky:Skoda_}). 

\hfill

$\bullet$ In the last section
(Section \ref{_stabi_final_Section_}), we use the results of 
 Section
\ref{_positi_hype_Section_}
(the positivity of a degree of a closed 2-form $\eta$ which is
a sum of a positive and an $SU(2)$-invariant form) 
to prove our main result. Given a sheaf $F$ admitting 
a connection with $SU(2)$-invariant curvature, 
we show that $F$ is a direct sum of stable sheaves.
This is done in the same way as one proves
that a Yang-Mills bundle is polystable.
We use the standard inequality between the
curvature of a bundle and a sub-bundle, which is
proven via the second fundamental form of
a sub-bundle (\cite{_Demailly_}, \cite{_Griffi_Harri_}).


\section{Hyperk\"ahler manifolds}
\label{_hyperka_Section_}


This Section contains a compression of 
the basic and best known results 
and definitions from hyperk\"ahler geometry, found, for instance, in
\cite{_Besse:Einst_Manifo_}, \cite{_Beauville_} and 
\cite{_Verbitsky:Hyperholo_bundles_}.

\hfill

\definition \label{_hyperkahler_manifold_Definition_} 
(\cite{_Besse:Einst_Manifo_}) A {\bf hyperk\"ahler manifold} is a
Riemannian manifold $M$ endowed with three complex structures $I$, $J$
and $K$, such that the following holds.
 
\begin{description}
\item[(i)]  the metric on $M$ is K\"ahler with respect to these complex 
structures and
 
\item[(ii)] $I$, $J$ and $K$, considered as  endomorphisms
of a real tangent bundle, satisfy the relation 
$I\circ J=-J\circ I = K$.
\end{description}

\hfill 

The notion of a hyperk\"ahler manifold was 
introduced by E. Calabi (\cite{_Calabi_}).

\hfill

Clearly, a hyperk\"ahler manifold has a natural action of
the quaternion algebra ${\Bbb H}$ in its real tangent bundle $TM$. 
Therefore its complex dimension is even.
For each quaternion $L\in \Bbb H$, $L^2=-1$,
the corresponding automorphism of $TM$ is an almost complex
structure. It is easy to check that this almost 
complex structure is integrable (\cite{_Besse:Einst_Manifo_}).

\hfill

\definition \label{_indu_comple_str_Definition_} 
Let $M$ be a hyperk\"ahler or hypercomplex
manifold, and $L$ a quaternion satisfying
$L^2=-1$. The corresponding complex structure 
on $M$ is called
{\bf an induced complex structure}. The $M$, considered as a K\"ahler
manifold, is denoted by $(M, L)$. In this case,
the hyperk\"ahler 
structure is called {\bf compatible
with the complex structure $L$}.

\hfill

\definition \label{_holomorphi_symple_Definition_} 
Let $M$ be a complex manifold and $\Theta$ a closed 
holomorphic 2-form over $M$ such that 
$\Theta^n=\Theta\wedge\Theta\wedge...$, is
a nowhere degenerate section of a canonical class of $M$
($2n=dim_\C(M)$).
Then $M$ is called {\bf holomorphically 
symplectic}.

\hfill

Let $M$ be a hyperk\"ahler manifold; denote the
Riemannian form on $M$ by $\langle \cdot,\cdot\rangle $.
Let the form $\omega_I := \langle I(\cdot),\cdot\rangle$ be the usual K\"ahler
form  which is closed and parallel
(with respect to the Levi-Civita connection). Analogously defined 
forms $\omega_J$ and $\omega_K$ are
also closed and parallel. 
 
A simple linear algebraic
consideration (\cite{_Besse:Einst_Manifo_}) shows that the form
\begin{equation}\label{_holo_symple_expli_Equation_}
\Theta:=\omega_J+\sqrt{-1}\omega_K
\end{equation}
is of
type $(2,0)$ and, being closed, this form is also holomorphic.
Also, the form $\Theta$ is nowhere degenerate, as another linear 
algebraic argument shows.
It is called {\bf the canonical holomorphic symplectic form
of a manifold M}. Thus, for each hyperk\"ahler manifold $M$,
and an induced complex structure $L$, the underlying complex manifold
$(M,L)$ is holomorphically symplectic. The converse assertion
is also true:

\hfill

\theorem \label{_symplectic_=>_hyperkahler_Proposition_}
(\cite{_Beauville_}, \cite{_Besse:Einst_Manifo_})
Let $M$ be a compact holomorphically
symplectic K\"ahler manifold with the 
holomorphic symplectic form
$\Theta$, a K\"ahler class 
$[\omega]\in H^{1,1}(M)$ and a complex structure $I$. 
Let $n=\dim_\C M$. Assume that
$\int_M \omega^n = \int_M (Re \Theta)^n$.
Then there is a unique hyperk\"ahler 
structure $(I,J,K,(\cdot,\cdot))$
over $M$ such that the cohomology class of the symplectic form
$\omega_I=(\cdot,I\cdot)$ is equal to $[\omega]$ and the
canonical symplectic form $\omega_J+\1\omega_K$ is
equal to $\Theta$.

\hfill

\ref{_symplectic_=>_hyperkahler_Proposition_} 
follows from the conjecture of Calabi, pro\-ven by
S.-T. Yau (\cite{_Yau:Calabi-Yau_}). 
\endproof

\hfill

Let $M$ be a hyperk\"ahler manifold. We identify the group $SU(2)$
with the group of unitary quaternions. This gives a canonical 
action of $SU(2)$ on the tangent bundle, and all its tensor
powers. In particular, we obtain a natural action of $SU(2)$
on the bundle of differential forms. 

\hfill

The following lemma is clear.

\hfill

\lemma \label{_SU(2)_commu_Laplace_Lemma_}
The action of $SU(2)$ on differential forms commutes
with the Laplacian.
 
{\bf Proof:} This is Proposition 1.1
of \cite{_Verbitsky:Symplectic_II_}. \endproof
 
\hfill

Thus, for compact $M$, we may speak of the natural action of
$SU(2)$ in cohomology.

\hfill

Further in this article, we use the following statement.

\hfill

\lemma \label{_SU(2)_inva_type_p,p_Lemma_} 
Let $\omega$ be a differential form over
a hyperk\"ahler or hypercomplex 
manifold $M$. The form $\omega$ is $SU(2)$-invariant
if and only if it is of Hodge type $(p,p)$ with respect to all 
induced complex structures on $M$.

\hfill

{\bf Proof:} Let $I$ be an induced complex structure,
and $\rho_I:\; U(1) \arrow SU(2)$ the corresponding
embedding, induced by the map 
$\R= \goth{u}(1) \arrow \goth{su}(2)$, $1\arrow I$.
The Hodge decomposition on $\Lambda^*(M)$ 
coincides with the weight decomposition of 
the $U(1)$-action $\rho_I$.
An $SU(2)$-invariant form is
also invariant with respect to $\rho_I$,
and therefore, has Hodge type $(p,p)$. Conversely,
if a $\eta$ is invariant with respect to $\rho_I$,
for all induced complex structures $I$, then
$\eta$ is invariant with respect to the
Lie group $G$ generated by these 
$U(1)$-subgroups of $SU(2)$.
A trivial linear-algebraic 
argument ensures that $G$ is the whole $SU(2)$.
This proves \ref{_SU(2)_inva_type_p,p_Lemma_}.
\endproof


\section{Hyperk\"ahler manifolds and stable bundles}
\label{_hyperho_Section_}


\subsection{Hyperholomorphic 
connections}

 Let $B$ be a holomorphic vector bundle over a complex
manifold $M$, $\nabla$ a  connection 
in $B$ and $\Theta\in\Lambda^2(M)\otimes End(B)$ be its curvature. 
This connection
is called {\bf compatible with a holomorphic structure} if
$\nabla_X(\zeta)=0$ for any holomorphic section $\zeta$ and
any antiholomorphic tangent vector field $X\in T^{0,1}(M)$. 
If there exists a holomorphic structure compatible with the given
Hermitian connection then this connection is called
{\bf integrable}.

\hfill
 
One can define a {\bf Hodge decomposition} in the space of differential
forms with coefficients in any complex bundle, in particular,
$End(B)$.

\hfill

\theorem \label{_Newle_Nie_for_bu_Theorem_}
Let $\nabla$ be a Hermitian connection in a complex vector
bundle $B$ over a complex manifold. Then $\nabla$ is integrable
if and only if $\Theta\in\Lambda^{1,1}(M, \End(B))$, where
$\Lambda^{1,1}(M, \End(B))$ denotes the forms of Hodge
type (1,1). Also,
the holomorphic structure compatible with $\nabla$ is unique.

{\bf Proof:} This is Proposition 4.17 of \cite{_Kobayashi_}, 
Chapter I.
$\blacksquare$

\hfill

This proposition is a version of Newlander-Nirenberg theorem.
For vector bundles, it was proven by M. Atiyah and R. Bott.

\hfill

\definition \label{_hyperho_conne_Definition_}
\cite{_Verbitsky:Hyperholo_bundles_}
Let $B$ be a Hermitian vector bundle with
a connection $\nabla$ over a hyperk\"ahler manifold
$M$. Then $\nabla$ is called {\bf hyperholomorphic} if 
the curvature of $\nabla$ is $SU(2)$-invariant.

\hfill

\example \label{_tangent_hyperholo_Example_} 
(Examples of hyperholomorphic bundles)

\begin{description}

\item[(i)]
Let $M$ be a hyperk\"ahler manifold, and $TM$ be its tangent bundle
equi\-p\-ped with the Levi--Civita connection $\nabla$. Consider a complex
structure on $TM$ induced from the quaternion action. Then $\nabla$
is a Hermitian connection
which is integrable with respect to each induced complex structure,
and hence, is hyperholomorphic.

\item[(ii)] For $B$ a hyperholomorphic bundle, all its tensor powers
are hyperholomorphic.

\item[(iii)] Thus, the bundles of differential forms on a hyperk\"ahler
manifold are also hyperholomorphic.

\end{description}


\subsection{Hyperholomorphic bundles 
and Yang-Mills connections.}
\label{_sta_bu_and_YM_Subsection_}


\definition\label{_degree,slope_destabilising_Definition_} 
Let $F$ be a coherent sheaf over
an $n$-dimensional compact K\"ahler manifold $M$. We define
{\bf the degree} $\deg(F)$ (sometimes the degree
is also denoted by $\deg c_1(F)$) as
\[ 
   \deg(F)=\int_M\frac{ c_1(F)\wedge\omega^{n-1}}{vol(M)}
\] 
and $\text{slope}(F)$ as
\[ 
   \text{slope}(F)=\frac{1}{\text{rank}(F)}\cdot \deg(F). 
\]
The number $\text{slope}(F)$ depends only on a
cohomology class of $c_1(F)$. 

Let $F$ be a torsion-free coherent sheaf on $M$
and $F'\subset F$ its proper subsheaf. Then $F'$ is 
called {\bf destabilizing subsheaf} 
if $\text{slope}(F') \geq \text{slope}(F)$

A coherent sheaf $F$ is called {\bf 
 stable}
\footnote{In the sense of Mumford-Takemoto}
if it has no destabilizing subsheaves. 
A coherent sheaf $F$ is called {\bf 
polystable} if it is a direct sum of stable sheaves of the same slope.
 
\hfill

Let $M$ be a K\"ahler manifold with a K\"ahler form $\omega$.
For differential forms with coefficients in any vector bundle
there is a Hodge operator $L: \eta\arrow\omega\wedge\eta$.
There is also a fiberwise-adjoint Hodge operator $\Lambda$
(see \cite{_Griffi_Harri_}).
 
\hfill

\definition \label{Yang-Mills_Definition_}
Let $B$ be a holomorphic bundle over a K\"ahler manifold $M$
with a holomorphic Hermitian connection $\nabla$ and a 
curvature $\Theta\in\Lambda^{1,1}\otimes End(B)$.
The Hermitian metric on $B$ and the connection $\nabla$
defined by this metric are called {\bf 
Yang-Mills} if 

\[
   \Lambda(\Theta)=constant\cdot \Id\restrict{B},
\]
where $\Lambda$ is a Hodge operator and $\Id\restrict{B}$ is 
the identity endomorphism which is a section of $End(B)$.

\hfill
 
A holomorphic bundle is called  {\bf indecomposable} 
if it cannot be decomposed onto a direct sum
of two or more holomorphic bundles.

\hfill

The following fundamental 
theorem provides examples of 
Yang-\-Mills \linebreak bundles.

\theorem \label{_UY_Theorem_} 
(Uhlenbeck-Yau)
Let B be an indecomposable
holomorphic bundle over a compact K\"ahler manifold. Then $B$ admits
a Hermitian 
Yang-Mills connection if and only if it is 
stable. Moreover, the Yang-Mills 
connection is unique, if it exists.
 
{\bf Proof:} \cite{_Uhle_Yau_}. \endproof

\hfill

\proposition \label{_hyperholo_Yang--Mills_Proposition_}
Let $M$ be a hyperk\"ahler manifold, $L$
an induced complex structure and $B$ be a complex vector
bundle over $(M,L)$. 
Then every 
hyperholomorphic connection $\nabla$ in $B$
is Yang-Mills and satisfies $\Lambda(\Theta)=0$
where $\Theta$ is a curvature of $\nabla$.
 
\hfill

{\bf Proof:} We use the definition of a hyperholomorphic 
connection as one with $SU(2)$-invariant curvature. 
Then \ref{_hyperholo_Yang--Mills_Proposition_}
follows from the

\hfill

\lemma \label{_Lambda_of_inva_forms_zero_Lemma_}
Let $\Theta\in \Lambda^2(M)$ be a $SU(2)$-invariant 
differential 2-form on $M$. Then
$\Lambda_L(\Theta)=0$ for each induced complex structure
$L$.\footnote{By $\Lambda_L$ we understand the Hodge operator 
$\Lambda$ adjoint to the multiplication by the K\"ahler form
associated with the complex structure $L$.}

{\bf Proof:} This is Lemma 2.1 of \cite{_Verbitsky:Hyperholo_bundles_}.
\endproof
 
\hfill

Let $M$ be a compact hyperk\"ahler manifold, $I$ an induced 
complex structure. 
For any 
stable holomorphic bundle on $(M, I)$ there exists a unique
Hermitian  Yang-Mills connection 
which, for some bundles, turns out to be hyperholomorphic. 
It is possible to tell exactly when
this happens.

\hfill

\theorem \label{_inva_then_hyperho_Theorem_}
Let $B$ be a 
stable holomorphic bundle over
$(M,I)$, where $M$ is a hyperk\"ahler manifold and $I$
is an induced complex structure over $M$. Then 
$B$ admits a compatible 
hyperholomorphic connection if and only
if the first two Chern classes $c_1(B)$ and $c_2(B)$ are 
$SU(2)$-invariant.\footnote{We use \ref{_SU(2)_commu_Laplace_Lemma_}
to speak of action of $SU(2)$ in cohomology of $M$.}

{\bf Proof:} This is Theorem 2.5 of
 \cite{_Verbitsky:Hyperholo_bundles_}. \endproof


\section{Hyperholomorphic sheaves}
\label{_hyperho_she_Section_}


In \cite{_Bando_Siu_}, S.  Bando and Y.-T. Siu developed 
machinery allowing one to apply the methods of Yang-Mills
theory to torsion-free coherent sheaves. In
\cite{_V:Hyperholo_sheaves_}, 
 their work was applied to generalise the results of
\cite{_Verbitsky:Hyperholo_bundles_}  
(see Section \ref{_hyperho_Section_})
to coherent sheaves. The first two subsections 
of this Section are a 
compilation of the results
and definitions of 
\cite{_Bando_Siu_} and \cite{_V:Hyperholo_sheaves_}.

\hfill

\subsection{Stable sheaves and Yang-Mills connections}
\label{_stable_refle_YM_Subsection_}

In this subsection, we repeat the basic definitions
and results from \cite{_Bando_Siu_} and \cite{_OSS_}.

\hfill

\definition\label{_refle_Definition_}
Let $X$ be a complex manifold, and $F$ a coherent sheaf on $X$.
Consider the sheaf $F^*:= \c Hom_{\calo_X}(F, \calo_X)$.
There is a natural functorial map 
$\rho_F:\; F \arrow F^{**}$. The sheaf $F^{**}$
is called {\bf a reflexive hull}, or {\bf 
reflexization}
of $F$. The sheaf $F$ is called {\bf reflexive} if the map
$\rho_F:\; F \arrow F^{**}$ is an isomorphism. 

\hfill

\remark\label{_refle_obtained_Remark_}
For all coherent sheaves $F$, the map
$\rho_{F^*}:\; F^* \arrow F^{***}$ is an isomorphism
(\cite{_OSS_}, Ch. II, the proof of Lemma 1.1.12).
Therefore, a 
reflexive hull of a sheaf is always 
reflexive. Moreover, a reflexive hull can be obtained by 
restricting to a non-singular set of $F$ subset  and taking the
pushforward (\cite{_OSS_}, Ch. II, Lemma 1.1.12). 
More generally, a reflexive sheaf is isomorphic
to a pushforward of its restriction
to an open set $M\backslash Z$, for
any complex analytic subset $Z\subset M$
of codimension at least 2.

\hfill

\lemma\label{_refle_pushfor_Lemma_}
Let $X$ be a complex manifold, $F$ a coherent sheaf on $X$,
$Z$ a closed analytic subvariety, $\codim Z\geq 2$, and
$j:\; (X\backslash Z) \hookrightarrow X$ the natural
embedding. Assume that the pullback $j^* F$ is
reflexive on $(X\backslash Z)$. Then the pushforward
$j_* j^* F$ is also reflexive.

\hfill

{\bf Proof:} This is \cite{_OSS_}, Ch. II, Lemma 1.1.12.
\endproof

\hfill

\lemma\label{_singu_codim_3_refle_Lemma_}
Let $F$ be a reflexive sheaf on $M$, and $X$ its singular set.
Then $\codim_M X\geq 3$

\hfill

{\bf Proof:} This is \cite{_OSS_}, Ch. II, 1.1.10. \endproof

\hfill

\claim 
Let $X$ be a K\"ahler manifold, and $F$ a torsion-free coherent sheaf over 
$X$. Then $F$ (semi)stable if and only if $F^{**}$ 
is (semi)stable.

{\bf Proof:} This is
\cite{_OSS_}, Ch. II, Lemma 1.2.4.
\endproof

\hfill

The admissible 
Hermitian metrics, introduced by  Bando and  Siu
in \cite{_Bando_Siu_}, play the role of the
ordinary Hermitian metrics for vector bundles.

\hfill

Let $X$ be a K\"ahler manifold.
In Hodge theory, one considers the operator 
$\Lambda:\; \Lambda^{p, q}(X) \arrow\Lambda^{p-1, q-1}(X)$
acting on differential forms on $X$, which is adjoint to the
multiplication by the K\"ahler form. This operator is also defined
on differential forms with coefficients in a bundle.
Consider a curvature $\Theta$ of a bundle $B$
as a 2-form with coefficients in $\End(B)$. Then
$\Lambda\Theta$ is a section of
$\End(B)$.

\hfill

\definition \label{_admi_metri_Definition_}
Let $X$ be a K\"ahler manifold, and $F$ a 
reflexive 
coherent sheaf over $X$. Let $U\subset X$ be the set of all
points at which $F$ is locally trivial. By definition,
the restriction $F\restrict U$ of $F$ to $U$ is a bundle.
An {\bf admissible metric} on $F$ is a Hermitian metric $h$
on the bundle $F\restrict U$ which satisfies the following
assumptions
\begin{description}
\item[(i)] the curvature $\Theta$ of $(F, h)$ is square integrable, and 
\item[(ii)] the corresponding section $\Lambda \Theta\in \End(F\restrict U)$ 
is uniformly bounded. 
\end{description}

\hfill

\definition \label{_Yang-Mills_sheaves_Definition_}
Let $X$ be a K\"ahler manifold, $F$ a 
reflexive sheaf over $X$, and $h$ an admissible metric on $F$.
Consider the corresponding Hermitian connection
$\nabla$ on $F\restrict U$. The metric $h$ and
the Hermitian connection $\nabla$ are called {\bf Yang-Mills}
if its curvature satisfies
\[ \Lambda \Theta\in \End(F\restrict U) = c\cdot \id 
\]
where $c$ is a constant and $\id$ the unit section 
$\id \in \End(F\restrict U)$. 

\hfill

One of the main results of \cite{_Bando_Siu_}
is the following analogue of the Uh\-len\-beck-\-Yau theorem
(\ref{_UY_Theorem_}).

\hfill

\theorem\label{_UY_for_shea_Theorem_}
Let $M$ be a compact K\"ahler manifold,  and $F$ a coherent
sheaf without torsion. Then $F$ admits an admissible 
 Yang--Mills 
metric if and only if $F$ is polystable. Moreover, if $F$
is stable, then this metric is unique, up to a constant
multiplier.

{\bf Proof:} \cite{_Bando_Siu_}, Theorem 3.
\endproof

\hfill

\remark
Clearly, the connection, corresponding to a metric on $F$,
does not change when the metric is multiplied by a scalar.
The  Yang--Mills metric on a polystable sheaf is unique up to
a component-wise multiplication by scalar multipliers.
Thus, the Yang--Mills connection of \ref{_UY_for_shea_Theorem_}
is unique.

\subsection{Stable hyperholomorphic sheaves over
hyperk\"ahler manifolds}
\label{_hyperholo_shea_def_Subsection_}

Let $M$ be a compact hyperk\"ahler manifold, $I$ an induced
complex structure, $F$ a torsion-free coherent sheaf
over $(M,I)$ and $F^{**}$ its 
reflexization. Recall that the cohomology of
$M$ are equipped with a natural $SU(2)$-action
(\ref{_SU(2)_commu_Laplace_Lemma_}). The motivation for
the following definition is \ref{_inva_then_hyperho_Theorem_} 
and \ref{_UY_for_shea_Theorem_}.

\hfill

\definition \label{_hyperho_shea_Definition_}
Assume that the first two Chern classes
of the sheaves $F$, $F^{**}$ are $SU(2)$-invariant.
Then $F$ is called {\bf stable hyperholomorphic} if
$F$ is stable. 
If $F$ is a direct sum of stable hyperholomorphic
sheaves, $F$ is called {\bf polystable hyperholomorphic}, 

\hfill

\remark \label{_slope_hyperho_Remark_}
The slope of a hyperholomorphic sheaf is zero, because
a degree of an $SU(2)$-invariant second cohomology class
is zero (\ref{_Lambda_of_inva_forms_zero_Lemma_}). 

\hfill

Let $M$ be a hyperk\"ahler manifold, $I$ an induced complex
structure, and $F$ a torsion-free sheaf over $(M,I)$.
Consider the natural $SU(2)$-action in the bundle
$\Lambda^i (M,B)$ of the differential $i$-forms
with coefficients in a vector bundle $B$. Let
$\Lambda^i_{inv}(M, B)\subset \Lambda^i (M, B)$
be the bundle of $SU(2)$-invariant $i$-forms.

\hfill

\definition \label{_hyperholo_co_Definition_}
Let $Z\subset (M, I)$ be a complex subvariety of
codimension at least 2, and $F$ a reflexive sheaf on $(M,I)$,
such that $F\restrict{M\backslash Z}$
is a bundle. Consider an admissible metric $h$ on 
$F\restrict{M\backslash Z}$, and let $\nabla$ be the associated
connection. Then $\nabla$ is called {\bf admissible hyperholomorphic} if its
curvature 
\[ \Theta_\nabla = 
   \nabla^2 \in 
   \Lambda^2\left(M, \End\left(F\restrict{M\backslash Z}\right)\right)
\]
is $SU(2)$-invariant, i. e. belongs to 
$\Lambda^2_{inv}\left(M, 
\End\left(F\restrict{M\backslash Z}\right)\right)$.

\hfill

\remark
This is the same definition as
\ref{_admissi_hy_Definition_}.

\hfill

\theorem\label{_hyperho_conne_exi_Theorem_}
Let $M$ be a compact hyperk\"ahler manifold, $I$
an induced complex structure and $F$ a reflexive sheaf on
$(M,I)$. Then $F$ admits a hyperholomorphic connection if
and only if $F$ is polystable hyperholomorphic.

\hfill

{\bf Proof:} This is \cite{_V:Hyperholo_sheaves_}, 
Theorem 3.19. \endproof


\subsection{Weakly hyperholomorphic sheaves}
\label{_weak_hype_Subsection_}

\definition \label{_wea_hh_Definition_}
Let $M$ be a hyperk\"ahler manifold, $I$ an induced complex
structure, and $F$ a torsion-free coherent sheaf on $M$.
Assume that outside of a closed complex analytic set $Z\subset (M,I)$,
$\codim _\C Z\geq 3$, the sheaf $F$ is smooth and
equipped with a connection $\nabla$. Assume, moreover, that
the curvature of $\nabla$ is $SU(2)$-invariant. 
Then $F$ is called {\bf weakly hyperholomorphic}.

\hfill

The main result of this paper is the following theorem.

\hfill

\theorem \label{_wea_hh_hh_Theorem_}
Let $M$ be a compact hyperk\"ahler manifold, $I$ an induced complex
structure, and $F$ a reflexive sheaf on $(M,I)$.
Assume that $F$ is weakly hyperholomorphic. Then
$F$ is polystable.

{\bf Proof:} See Section \ref{_stabi_final_Section_}.
 \endproof

\hfill

\remark
From \ref{_wea_hh_hh_Theorem_} it follows that
all stable summands $F_i$ of $F$ are hyperholomorphic 
(see \ref{_SU(2)-inv_c_i_Remark_}).

\hfill

By \ref{_hyperho_conne_exi_Theorem_},
$F$ admits a unique admissible hyperholomorphic
Yang-Mills connection $\nabla_1$. However,
we do not know whether $\nabla=\nabla_1$ or not.

\hfill

\conjecture \label{_wea_hh_conn_hh_Conjecture_}
Under assumptions of \ref{_wea_hh_hh_Theorem_},
the connection $\nabla$ is admissible.

\hfill

Clearly, \ref{_wea_hh_hh_Theorem_} is implied by 
\ref{_hyperho_conne_exi_Theorem_} and \ref{_wea_hh_conn_hh_Conjecture_}.

\hfill

\example 
Let $B$ be a hyperholomorphic bundle on a product
$M_1\times M_2$ of two hyperkaehler manifolds,
and $M_1 \times M_2 \stackrel \pi \arrow M_1$ the
projection map. From the usual twistor argument 
(see e. g. \cite{_NHYM_} or \cite{_BBR_3}) it follows 
that a derived direct image $R^i\pi_* B$ admits
a hyperholomorphic connection outside of its
singularities. If, in addition, $M_1$ is generic
in its deformation class, all its subvarieties
have even codimension. In particular, $R^i\pi_* B$
is smooth outside of codimension 2, and
its reflexization $(R^i\pi_* B)^{**}$
is weakly hyperholomorphic. The stability
of direct images of stable bundles is
an important question which is partially
solved by \ref{_wea_hh_hh_Theorem_}.


\section{Positive forms and hyperholomorphic connections}
\label{_posi_c_2_Section_}


To justify \ref{_wea_hh_conn_hh_Conjecture_},
we prove it for sheaves with isolated singularities. 

\subsection{Singularities of positive closed forms}

\definition
Let $M$ be a complex manifold, and $\eta$ a real-valued 
$(p,p)$-form on $M$. Then $\eta$ is called {\bf positive}
if for any $p$-tuple of vector fields 
\[ \alpha_1, ..., \alpha_1\in \Lambda^{1,0}(M),\]
we have
\[ (\1)^p \eta(\alpha_1, \bar \alpha_1, \alpha_2, 
\bar\alpha_2, ...)\geq 0.
\]
For an excellent exposition of the theory of 
positive forms and currents, see \cite{_Demailly_}.
Further on, we shall need the following important lemma.

\hfill

\lemma\label{_Sibony_Lemma_}
Let $M$ be a K\"ahler manifold, $Z\subset M$ a closed complex
subvariety, $\codim Z>p$, and $\eta$ a closed positive $(p,p)$-form
on $M\backslash Z$. Then $\eta$ is locally $L^1$-integrable.

\hfill

{\bf Proof:} \cite{_Sibony_}. \endproof

\subsection{Weakly holomorphic sheaves with isolated singularities}

The following proposition is not used anywhere in this paper.
We include it to justify \ref{_wea_hh_conn_hh_Conjecture_},
and, ultimately - to support \ref{_wea_hh_hh_Theorem_}
with a simple and convincing argument, albeit
valid only in a special case.

\hfill

\proposition \label{_iso_sing_curv_L^2_Proposition_}
Let $M$ be a compact hyperk\"ahler manifold, $I$ an induced complex
structure, $F$ a torsion-free coherent sheaf on $(M,I)$ with isolated
singularities, and $Z\subset (M,I)$ a finite set 
containing the singular points of $F$. Consider
a hyperholomorphic connection $\nabla$ on
$F\restrict {M\backslash Z}$.
Then $\nabla$ is admissible, in the sense of
\ref{_admi_metri_Definition_}.\footnote{Since $\nabla$ is Yang-Mills,
it is admissible if and only if its curvature is square-integrable.}

\hfill

{\bf Proof:} Consider the curvature as a form 
$\Theta \in \Lambda^{1,1}(M)\otimes \goth{su}(B)$,
where $ \goth{su}(B)$ denotes the Lie algebra
of traceless skew-Hermitian endomorphisms of $B$
Let $r_2\in \Lambda^{2,2}(M)$ be the form 
\[ r_2:=
   \frac{\1}{{2\pi}^2}\Tr(\Theta\wedge \Theta)
\]
representing, by Gauss-Bonnet formula, 
the cohomology class $2 c_2(B) - \frac{n-1}{n} c_1(B)^2$.
Let $\omega$
be the K\"ahler form of $(M, I)$, and $\Vol(M)$
is volume form.
By the Hodge-Riemann relations, for $\Lambda\Theta=0$, we have
\[ 
  \Tr(\Theta\wedge \Theta) \wedge \omega^{n-2} = 
  (4 ( n^2-n))^{-1} \| \Theta \|^2\Vol(M)
\]
(see, e.g. \cite{_Bando_Siu_} or \cite{_V:Hyperholo_sheaves_}).
Therefore,  $\Theta$ is square-integrable if and only if
$r_2 \wedge \omega^{n-1}$ lies in $L^1(M)$.

\hfill

We use the following fundamental lemma; 
we shall complete the proof of 
\ref{_iso_sing_curv_L^2_Proposition_} at the end of
this section.

\hfill

\lemma \label{_c_2_wedge_omega^n-3_positive_Lemma_}
Let $M$, $\dim_\C M =n$, $n> 2$ 
be a hyperk\"ahler manifold, $(B, \nabla)$ a 
holomorphic bundle with a hyperholomorphic connection, 
$\Theta$ its curvature, and $r_2:= \frac{\1}{{2\pi}^2}\Tr(\Theta\wedge \Theta)$
be the corresponding 4-form, representing (by Gauss-Bonnet)
$2 c_2(B) - \frac{n-1}{n} c_1(B)^2$.
Consider the $(n-1, n-1)$-form $r_2 \wedge \omega^{n-3}$.
Then $r_2 \wedge \omega^{n-3}$ 
is positive.\footnote{Since the statement
of \ref{_c_2_wedge_omega^n-3_positive_Lemma_} is local,
it is also true for an admissible hyperholomorphic
connection on a reflexive sheaf.}

\hfill

{\bf Proof:} The statement of \ref{_c_2_wedge_omega^n-3_positive_Lemma_}
is essentially linear-algebraic. Let $x\in M$ be a point, and
$z_1, ... z_n, \bar z_1, ... \bar z_n$ a standard 
basis in the complexified cotangent space $T_x^* M\otimes \C$,
such that the K\"ahler form is written as
\[ \omega\restrict x = \1\sum_i z_i \wedge \bar z_i. 
\]
Denote by $w_{ij}\in \Lambda^{n-1, n-1}(M)$ the form
\[ w_{ij}:= 
  z_1\wedge ... \check z_i ... \wedge z_n 
  \wedge \bar z_1\wedge ... \check {\bar z_j} ... \wedge \bar z_n 
\]
where $z_1\wedge ... \wedge\check z_i \wedge ... \wedge z_n$ denotes a product
of all $z_k$ {\em except} $z_i$.

Write $r_2 \wedge \omega^{n-3}$ in this basis as
\[ 
   r_2 \wedge \omega^{n-3} = \sum_{i,j} B_{ij} w_{ij},\ \  B_{ij}\in \R
\]
To prove
\ref{_c_2_wedge_omega^n-3_positive_Lemma_},
we need to show that 
\begin{equation}\label{_B_ii_posi_Equation_}
  B_{ii}\geq 0,\  i = 1, ... n.
\end{equation}
Indeed, positivity of a real $(n-1,n-1)$-form $\nu$
is equivalent to the positivity of a product
$\1\nu\wedge z\wedge\bar z$, for each $z\in \Lambda^{1,0}(M)$.
Choosing the basis $z_1, ..., z_n \in \Lambda^{1,0}(M)$
in such a way that $z_i=z$, we find that whenever $B_{1,1}\geq 0$
for each orthonormal basis $z_1, ..., z_n$, the form $r_2 \wedge \omega^{n-3}$
is positive. 

Write $\Theta\restrict x$ as
\[ \Theta\restrict x = \sum_{i,j} z_i \wedge\bar z_j A_{ij}
\]
where $A_{ij}\in \goth{su}\bigg(B\restrict x\bigg)$.
An easy calculation implies
\begin{equation}\label{B_ii_explicit_Equation_}
   B_{ii} = - (n-3)! \sum_{k, l} \Tr(A_{kl}^2) +
   (n-3)! \sum_{k, l} \Tr(A_{kk}A_{ll}) \ \ \  k\neq l\neq i
\end{equation}
(the sum is performed over all $k, l = 1, ... n$, satisfying
$k\neq l\neq i$).

The first summand of the right hand side of \eqref{B_ii_explicit_Equation_}
is non-negative, because $A_{kl}\in \goth{su}\bigg(B\restrict x\bigg)$,
and the Killing form on $\goth{su}\bigg(B\restrict x\bigg)$ is negative definite.

To prove \eqref{_B_ii_posi_Equation_}, it 
remains to show that the second summand of 
the right hand side of \eqref{B_ii_explicit_Equation_}
is non-negative:
\begin{equation}\label{_C_ii_defi_Equation_}
  C_{ii}:= \sum_{k, l}\Tr(A_{kk}A_{ll})\geq 0 \ \ \ (k\neq l\neq i).
\end{equation}
Clearly, 
\begin{equation} \label{_C_via_sum_A_kk_A_ii_i_fix_Equation_}
    C_{ii} = \sum_{k\neq l} \Tr(A_{kk}A_{ll}) - 
            2 \left(\sum_{k= 1, ..., i-1, i+1, ... n} \Tr(A_{kk}A_{ii})\right).
\end{equation}
Since $\Lambda(\Theta)=0$, we have
\begin{equation}\label{_sum_A_ii_zero_} 
  \sum_{k=1}^n A_{kk}=0. 
\end{equation}
Therefore, 
\begin{equation} \label{_sum_A_kk_A_ii_i_fix_Equation_}
  \sum_{k= 1, ..., i-1, i+1, ... n} \Tr(A_{kk}A_{ii}) = - \Tr (A_{ii}^2).
\end{equation}
Plugging \eqref{_sum_A_kk_A_ii_i_fix_Equation_} into 
\eqref{_C_via_sum_A_kk_A_ii_i_fix_Equation_}, we obtain
that \eqref{_sum_A_ii_zero_}  gives 
\begin{multline}\label{_C_ii_long_Equation_} 
  C_{ii} = \sum_{k\neq l} \Tr(A_{kk}A_{ll}) + 2 \Tr (A_{ii}^2) =
            -\sum_k \Tr (A_{kk}^2) + 2 \Tr (A_{ii}^2) \\
          =-\left(\sum_{k= 1, ..., i-1, i+1, ... n} \Tr (A_{kk}^2)\right) +  
           \Tr (A_{ii}^2)
\end{multline}

The formula \eqref{_C_ii_long_Equation_} was 
obtained using only the Yang-Mills property of the connection;
it is true for all K\"ahler manifolds. Now recall that
$\nabla$ is hyperholomorphic. 
Renumbering the basis $z_1, ... z_n$, we may assume 
that the number $i$ is odd. Fix the standard quaternion
triple $I, J, K$. This fixes a choice of
a holomorphic symplectic form (\ref{_holo_symple_expli_Equation_}). 
Changing $z_1, .. z_{i-1}, z_{i+1}, ... z_n$
if necessary, we may also assume that the 
holomorphic symplectic form is written as
\[ \Omega = z_1 \wedge z_2 + ... + z_i\wedge z_{i+1} \wedge ...
 + z_{n-1}\wedge z_{n}
\]
Let $J\in SU(2)$ be an operator of $SU(2)$ given by $J\in \Bbb H$. 
An easy calculation insures that $J$ maps the 2-form 
$z_i\wedge \bar z_{i}$ to $-z_i\wedge \bar z_{i}$.
Since $\Theta$ is $SU(2)$-invariant, we obtain that
$A_{ii} = - A_{i+1,i+1}$. Plugging this into 
\eqref{_C_ii_long_Equation_}, we find
\[
  C_{ii} = -\sum_{k\neq i, i+1}\Tr (A_{kk}^2)
\]
Since the Killing form is negative definite, the
number $C_{ii}$ is non-negative.
This proves \ref{_c_2_wedge_omega^n-3_positive_Lemma_}.
\endproof

\hfill

Return to the proof of \ref{_iso_sing_curv_L^2_Proposition_}.
We have shown that $r_2 \wedge \omega^{n-3}$
is a positive $(n-1, n-1)$-form; this form is also closed,
and smooth outside of a complex analytic subset
of codimension $n$.
By \ref{_Sibony_Lemma_}, such form is
$L^1$-integrable. We proved \ref{_iso_sing_curv_L^2_Proposition_}.
\endproof


\section{$SU(2)$-invariance of the Chern class}
\label{_SU(2)_inv_Section_}

\subsection{The Dolbeault spectral sequence and the Hodge filtration}

Further on in this section, we shall need the following proposition

\hfill

\proposition\label{_H^i_coh_codim_Proposition_}
Let $X$ be a complex manifold, and $Z\subset X$ a
real analytic subvariety admitting a stratification by smooth
real analytic subvarieties of even dimension.
Assume that $\codim _\R Z\geq 2m$, $m\geq 2$. Let 
$U= X\backslash Z$, and let $B$ be a holomorphic
bundle on $X$. Consider the natural map
of holomorphic cohomology $\phi:\; H^n(X, B)\arrow H^n(U, B)$.
Then $\phi$ is an isomorphism, for $n\leq m-2$.

\hfill

\remark
For $n=0$, \ref{_H^i_coh_codim_Proposition_} becomes the well
known Hartogs theorem. For $Z$ complex analytic,
\ref{_H^i_coh_codim_Proposition_} is also well known
(\cite{_Scheja_}). Further on, we shall use 
\ref{_H^i_coh_codim_Proposition_} when $X$ is
a hyperk\"ahler manifold with an induced complex structure
$J$, and $Z\subset M$ a complex analytic subvariety of
$(M,I)$.

\hfill

{\bf The proof of \ref{_H^i_coh_codim_Proposition_}}.

Using the Meyer-Vietoris exact sequence, we find that it suffices
to prove \ref{_H^i_coh_codim_Proposition_} when $X$ is an open ball
and the bundle $B$ is trivial. 
Using induction by $\dim Z$, we may also assume that
$Z$ is smooth (otherwise, prove \ref{_H^i_coh_codim_Proposition_}
for $X = X\backslash Sing(Z)$, and then apply
\ref{_H^i_coh_codim_Proposition_} to the pair
($X$, $Sing(Z)$).

Shrinking $X$ further and applying Meyer-Vietoris,
 we may assume that there exists a smooth holomorphic map $f:X \arrow Y$ 
inducing a real analytic isomorphism $f:\; Z\stackrel \sim \arrow Y$.
Shrinking $X$ again, we assume that $X= Y \times B$,
where $B$ is an open ball in $\C ^r$, $r\geq m$,
and $f:\; X \arrow Y$ is the standard projection map.
Denote by $g:\; X \arrow B$ the other standard projection map.

Consider the K\"unneth decomposition of the differential forms

\begin{equation}\label{_Kunneth_Equation_}
\Lambda^{0,a}(U) = \bigoplus_{b+c=a}
  f^*\Lambda^{0,b}(Y) \otimes_{C^\infty(X\backslash Z)} g^*\Lambda^{0,c}(B)
\end{equation}

Decomposing the Dolbeault complex of $U$ in accordance with
\eqref{_Kunneth_Equation_}, we obtain 
\begin{equation}\label{_Kunneth_diffe_Equation_}
\bar \6 = \bar \6_B + \bar\6_Y,
\end{equation}
where $\bar \6_B$ , $\bar\6_Y$ are the Dolbeault differentials on
$ \Lambda^{0,*}(B)$, $\Lambda^{0,*}(Y)$. Consider the bicomplex spectral
sequence, associated with the bicomplex \eqref{_Kunneth_Equation_}
and the decomposition \eqref{_Kunneth_diffe_Equation_}.
The cohomology of the complex $(f^* \Lambda^{0,*}(B),\bar\6_B)$ 
form a $C^\infty$-bundle on $Y$, 
with the fibers in $y\in Y$ identified with 
$H^{q}(\calo_{f^{-1}(y)})$.
Therefore, the 
$E_1$-term of this spectral sequence, $H^*(\Lambda^{0,*}(U),\bar \6_B)$
can be identified with the space of global sections of a 
graded $C^\infty$-bundle $R^{p,q}$ on $Y$, 
\[ R^{p,q}\restrict y = 
   H^{q}(\calo_{f^{-1}(y)})\otimes_\C \Lambda^{0,*}(T_y Y).
\]
By construction, the fibers $f^{-1}(y)$ are isomorphic to an
open ball without a point: $B^r\backslash {pt}$.
The cohomology of the structure sheaf on
$B^r\backslash {pt}$ are well known; in particular,
we have an isomorphism $H^i(\calo_{f^{-1}(y)}) \cong H^i(B)$,
$i\leq r-2$ (\cite{_Scheja_}).
Therefore, the natural functorial
morphism from cohomology of $X$ to the
cohomology of $U$ induces an isomorphism 
\begin{equation}\label{_spectra_equiv_Equation_}
E_1^{p,q}(X)\cong E_1^{p,q}(U),
\end{equation}
 for $q\leq r-2$.
This spectral sequence converges to $H^*(U, \calo_U)$.
Therefore, \eqref{_spectra_equiv_Equation_}
implies an isomorphism $H^n(\calo_X)\cong H^n(\calo_U)$,
for $n\leq m-2$. We have proved \ref{_H^i_coh_codim_Proposition_}.
\endproof

\hfill

\corollary\label{_Hodge_filt_Corollary_}
In assumptions of \ref{_H^i_coh_codim_Proposition_},
consider the natural map 
$H^{2n}(X) \stackrel \rho\arrow H^{2n}(U)$.
Since $n\leq m-2$, $\codim_\R  Z> 2n$, and
$\rho$ is an isomorphism. Let 
$F^0(X) \subset F^1(X) \subset ... \subset F^{2n}(X)=H^{2n}(X)$,
$F^0(U) \subset F^1(U) \subset ... \subset F^{2n}(U)=H^{2n}(U)$
be the Hodge filtration on the cohomology of $X$ and $U$.
Assume that $X$ is compact and K\"ahler. Then 
$\rho$ induces an isomorphism
\[ \rho:\; F^i(X) \arrow F^i(U)
\]
for $i\leq n$.

\hfill

{\bf Proof:} Consider the $E_2$-term of the Dolbeault spectral
sequence. Since $E_2^{p,q}(X)= H^q(\Omega^p(X))$ and
$E_2^{p,q}(U)= H^q(\Omega^p(U))$, the restriction map
$E_2^{2n-i,i}(X) \arrow E_2^{2n-i,i}(U)$ is an isomorphism
for $i\leq n$ (\ref{_H^i_coh_codim_Proposition_}).
By definition, $F^r$ is a union of all
elements of $\oplus_{i\leq r} E_2^{2n-i,i}$
which survive under the higher differentials,
up to the images of these differentials.
On the other hand, $X$ is K\"ahler and compact,
hence the Dolbeault spectral
sequence of $X$ degenerates in $E^2$. Therefore,
 the map $\rho:\; F^i(X) \arrow F^i(U)$
is surjective. It is injective because 
$H^{2n}(X) \stackrel \rho\arrow H^{2n}(U)$
is an isomorphism.

\endproof

\hfill

\corollary\label{_2-forms_corollary_}
Let $X$ be a compact K\"ahler manifold, $Z\subset X$ a
real analytic subvariety admitting a stratification by smooth
real analytic subvarieties of even dimension,
$\codim_\R Z\geq 6$, and $U:= X\backslash Z$.
Given a closed $(1,1)$-form
$\eta$ on $U$, the corresponding cohomology class
$[\eta]\in H^2(U)=H^2(X)$ has Hodge type $(1,1)$.

\hfill

{\bf Proof:} Since the $\eta$ is a $(1,1)$-form, 
the cohomology class $[\eta]\in H^2(U)$
belongs to the $F^1(U)$-term of the Hodge filtration.
By \ref{_Hodge_filt_Corollary_}, 
\begin{equation}\label{_[eta]_in_F^1_Equation_}
[\eta]\in F^1(X)= H^{2,0}(X) \oplus H^{1,1}(X).
\end{equation}
Replacing the complex structure $I$ on $X$ by $-I$,
 we obtain $H^{2,0}(X,I)= H^{0,2}(X, -I)$.
Applying the same argument to the 
cohomology class $[\eta]$ on $(X,-I)$,
we obtain
\begin{equation}\label{_[eta]_in_F^1(X,-I)_Equation_}
[\eta]\in F^1(X, -I)= H^{1,1}(X) \oplus H^{0,2}(X).
\end{equation}
Comparing 
\eqref{_[eta]_in_F^1_Equation_} and 
\eqref{_[eta]_in_F^1(X,-I)_Equation_},
we obtain \ref{_2-forms_corollary_}.
\endproof

\subsection{Closed $SU(2)$-invariant forms}

The main result of this section is the following theorem.

\hfill

\theorem\label{_wea_hyhe_c_1_SU(2)_inv_Theorem_}
Let $M$ be a compact hyperk\"ahler manifold, 
$I$ an induced complex structure and $F$ 
a weakly hyperholomorphic sheaf on $(M,I)$.
Then the Chern class $c_1(F)$ is $SU(2)$-invariant. 

\hfill

{\bf Proof:} Let $Z\subset (M,I)$ be the singular set of $F$.
By \ref{_singu_codim_3_refle_Lemma_}, $\codim _\C Z\geq 3$. 
Consider the form $\eta :=\Tr \Theta$ on $U:= M\backslash Z$,
where $\Theta$ is the curvature of $F\restrict U$.
Since $\codim _\C Z\geq 3$, we have $H^2(M)=H^2(U)$.
Clearly, the cohomology class $[\eta]\in H^2(U)=H^2(M)$
is equal to $c_1(F)$.

Let $L$ be an arbitrary induced complex structure on $M$.
Applying \ref{_2-forms_corollary_} to
the K\"ahler manifold $(M,L)$ and a closed $(1,1)$-form
$\eta$ on $M\backslash Z$, we find that
the cohomology class $[\eta]=c_1(F)$ is of
type $(1,1)$ with respect to $L$.
By \ref{_SU(2)_inva_type_p,p_Lemma_}, 
$c_1(F)$ is $SU(2)$-invariant. 
\endproof


\section{Positivity and hyperk\"ahler geometry}
\label{_positi_hype_Section_}


Let $M$ be a compact hyperk\"ahler manifold, $\dim_\R M = 2n$,
$I$ an induced complex structure, 
and $\omega_I$ the corresponding K\"ahler form
(Section \ref{_hyperka_Section_}). Given $\eta \in H^2(M)$, 
we define
\[ \deg_I(\eta) := \int_M \eta \wedge \omega_I^{n-1}.
\]
A 2-form is called {\bf pure of weight 2} if it is a sum of
forms which lie in 3-dimensional irreducible
$SU(2)$-subrepresentations of $\Lambda^2(M)$.

The aim of this section is the following theorem.

\hfill

\theorem \label{_degree_of_posi+inva_Theorem_}
In the above assumptions, let $Z\subset (M,I)$, $\codim_\C Z\geq 3$
be a closed complex subvariety, and 
$H^2(M) \arrow H^2(M\backslash Z)$ 
the induced isomorphism. Consider a closed (1,1)-form
$\eta$ on $(M\backslash Z, I)$, and let $[\eta]$ be its
cohomology class in $H^2(M\backslash Z) = H^2(M)$. Assume that
$\eta$ admits a decomposition $\eta = \eta_0 +\eta_+$,
where $\eta_0$ is $SU(2)$-invariant, and $\eta_+$ is positive
and pure of weight 2. Then $\deg_I ([\eta]) \geq 0$, and 
the equality is reached if and only if $\eta_+=0$.

\hfill

\remark 
The ``if'' part of the last statement is an immediate consequence
of \ref{_Lambda_of_inva_forms_zero_Lemma_}, \ref{_SU(2)_inva_type_p,p_Lemma_}  
and \ref{_2-forms_corollary_}.

\hfill

The proof of \ref{_degree_of_posi+inva_Theorem_}
takes the rest of this Section. 

\hfill

\lemma \label{_decompo_eta_w-K_SU(2)_Lemma_}
Let $M$ be a hyperk\"ahler manifold, $I, J, K$
the standard triple of induced complex structures,
and $\eta$ a $(1,1)$-form on $(M,I)$. Consider the decomposition
$\eta=\eta_+ + \eta_0$, where $\eta_0$ is $SU(2)$-invariant,
and $\eta_+$ is pure of weight 2. Consider the Hodge decomposition
\begin{equation}\label{_eta_wrt_K_Equation_}
\eta= \eta^{2,0}_K + \eta^{1,1}_K +\eta^{0,2}_K
\end{equation}
associated with $K$. Then
\begin{description}
\item[(i)]
\begin{equation}\label{_eta_wrt_K_and_SU(2)_Equation_}
\eta^{1,1}_K= \eta_0, \text{\ \ and \ \ } \eta^{2,0}_K  +\eta^{0,2}_K = \eta_+.
\end{equation}
\item[(ii)] 
The correspondence $\eta_+\arrow \eta^{2,0}_K$ induces an isomorphism 
between the bundle
$\Lambda^{2,0}_K(M)$ and the bundle $\Lambda^{1,1}_{I,+}(M)$
of (1,1)-forms of weight 2. 

\item[(iii)]
Moreover, this identification
maps the real structure $\eta_K^{2,0}\arrow I(\bar \eta_K^{2,0})$
to the standard real structure on $\Lambda^{1,1}_{I,+}(M)$.
\end{description}

{\bf Proof:} Let $K$ act on $\Lambda^*(M)$ 
multiplicatively as follows
\[ K(x_1\wedge x_2 \wedge ...) := K(x_1)\wedge K(x_2)\wedge ...
\]
Denote the action of $I$ in the same way.
Since the eigenvalues of $K$ on $\Lambda^1(M)$ are $\pm\1$,
the operator $K$ has eigenvalues $\pm 1$ on $\Lambda^2(M)$:
it acts on $\Lambda^{1,1}_K(M)$ as $1$, and on 
$\Lambda^{2,0}_K(M)\oplus \Lambda^{0,2}_K(M)$ 
as $-1$:
\begin{equation}\label{_K_eigenva_Equation_}
K\restrict {\Lambda^{1,1}_K(M)}=1, \ \ \ 
K\restrict{\Lambda^{2,0}_K(M)\oplus \Lambda^{0,2}_K(M)} =-1
\end{equation}
The central element $c\in SU(2)$ acts trivially
on $\Lambda^2(M)$; therefore, $I$ and $K$ commute on 
$\Lambda^2(M)$. We obtain that $K$ preserves the fixed
space of $I$:
\[ K(\Lambda^{1,1}_I(M)) = \Lambda^{1,1}_I(M).
\]
If a 2-form $\eta$ is fixed by $I$ and $K$,
it is also fixed by $K \circ I=J$. By
\eqref{_K_eigenva_Equation_}, this means that
$\eta$ is of type $(1,1)$ with respect to
$I$, $J$ and $K$. 
A simple linear-algebraic argument implies that $\eta$
is of type $(1,1)$ with respect to all induced complex structures.
By \ref{_SU(2)_inva_type_p,p_Lemma_},
this implies that $\eta$ is $SU(2)$-invariant.
We proved the first equation of 
\eqref{_eta_wrt_K_and_SU(2)_Equation_}. 

Now, if $\eta$ is pure of weight 2, 
$K$ acts on $\eta$ as on the K\"ahler form $\omega_I$;
it is easy to check, then, that $K(\eta)=-\eta$.
This implies that $\eta \in \Lambda^{2,0}_K(M)\oplus \Lambda^{0,2}_K(M)$.
Conversely, if $\eta$ belongs to $\Lambda^{2,0}_K(M)$ or
$\Lambda^{0,2}_K(M)$, it has weight $\pm 2$ with respect to the
Cartan algebra element, corresponding to $\1 K\in \goth{su}(2)\subset \H$;
therefore, $\eta$ is pure of weight 2. 
This proves \eqref{_eta_wrt_K_and_SU(2)_Equation_}.
This identifies the bundles
$\Lambda^{2,0}_K(M)$, $\Lambda^{0,2}_K(M)$
and the bundle $\Lambda^{1,1}_{I,+}(M)$ of $(1,1)$-forms
of weight 2. This identification is compatible
with the real structure on forms, and this gives the 
last assertion of \ref{_decompo_eta_w-K_SU(2)_Lemma_}.
\endproof

\hfill

\definition
Let $M$ be a hyperk\"ahler manifold, $I,J,K$ the standard 
triple of induced complex structures, and $\rho\in \Lambda^{2,0}_K(M)$
a (2,0)-form on $(M,K)$ satisfying $I(\rho)=\bar\rho$.
Consider the real part $\Re(\rho)$ of $\rho$.
Then 
\[ 2I(\Re(\rho)) = I(\rho) + I(\bar\rho)= \rho + \bar\rho
   = 2\Re(\rho).
\]
Therefore,
$\Re(\rho)$ lies inside $\Lambda^{1,1}_I(M, \R)$.
We say that $\rho$ is {\bf $K$-positive} if the form
$Re(\rho)$ is positive on $(M,I)$.

\hfill

\ref{_degree_of_posi+inva_Theorem_} is an immediate
corollary of the following Proposition, which is
a hyperk\"ahler version of \ref{_Sibony_Lemma_}.

\hfill

\proposition\label{_closed_K-positive_is_L^1_Proposition_}
Let $M$ be a hyperk\"ahler manifold, $I,J,K$ the standard 
triple of induced complex structures, $Z\subset (M,I)$
a compact complex subvariety, $\codim_\C Z\geq 3$, and 
$\rho\in \Lambda^{2,0}_K(M)$ a $\6_K$-closed 
$(2,0)$-form on $(M\backslash Z,K)$, which
satisfies $\rho= I(\bar\rho)$.
Assume that $\rho$ is $K$-positive. 
Then (i) $\rho$ is locally $L^1$-integrable on $M$,
and (ii) $\6_K$-closed as a current on $M$.

\hfill

We prove \ref{_closed_K-positive_is_L^1_Proposition_}
in Section \ref{_integra_K-posit_Section_}.
Let us show how to deduce 
\ref{_degree_of_posi+inva_Theorem_} from 
\ref{_closed_K-positive_is_L^1_Proposition_}.

\hfill

Let $\Omega_K:= \omega_I + \1 \omega_J$ be the holomorphic symplectic
form on $(M,K)$, and $n:= \dim_\H M$.
Consider the $(4n-2)$-form 
\[ E:= \Omega_K^{n-1} \wedge \bar \Omega_K^n.
\]
Further on, we shall need the following lemma.

\hfill

\lemma\label{_E_and_degree_Lemma_}
Let $[\eta]\in H^{1,1}_I(M)$
be a cohomology class, on a compact hyperk\"ahler manifold
$(M,I,J,K)$, $\dim_\R M = 4n$.
Then 
\begin{equation}\label{_degree_via_E_Equaion_}
\lambda_n\deg_I [\eta] =\Re\left( \int_M [\eta]_K^{2,0}\wedge E\right),
\end{equation}
where $\lambda_n$ is a positive constant, depending on $n$ only,
$[\eta]_K^{2,0}$ is a $(2,0)$-part of $[\eta]$
with respect to $K$, and $E$ the $(4n-2)$-form constructed
above.

\hfill

{\bf Proof:} By construction, $[\eta]_K^{2,0}\wedge E=[\eta] \wedge E$.
Also, $\deg_I [\eta]=\int_M[\eta]\wedge \omega_I^{2n-1}$.
Since $[\eta]\in H^{1,1}_I(M)$, we have
\[ [\eta] \wedge E = [\eta] \wedge E^{2n-1,2n-1}_I,
\]
where $E^{2n-1,2n-1}_I$ is $(2n-1, 2n-1)$-part of $E$,
taken with respect to $I$. Then 
\eqref{_degree_via_E_Equaion_} is implied by the relation
\begin{equation}\label{_E_via_omega_Equation_} 
E^{2n-1,2n-1}_I= \lambda_n^{-1} \omega_I^{2n-1}
\end{equation}
This relation is proven by direct calculation
(see e.g. \cite[Section 3]{_Verbitsky:Symplectic_I_}, where
the whole algebra generated by $\omega_I$, $\omega_J$ and $\omega_K$
is explicitly calculated). 
\endproof

\hfill

Given a $\6_K$-exact 2-form $\mu=\6_K\mu'$,
we have
\[ d(E\wedge \mu') = E\wedge \mu.
\]
Therefore, the number
\[ E(\rho) := \int_M \rho \wedge E
\]
depends only on the $\6_K$-Dolbeault 
cohomology class of $\rho \in \Lambda^{2,0}(M)$.

\hfill

Return to the assumptions of \ref{_degree_of_posi+inva_Theorem_}.
The form $\eta_K^{2,0}$ by construction
satisfies assumptions of
\ref{_closed_K-positive_is_L^1_Proposition_}. 
Therefore,  $\eta_K^{2,0}$
 is locally $L^1$-integrable. We shall interpret
$\deg_I [\eta]$ as integral 
$\lambda_n^{-1}\int_M \eta_K^{2,0}\wedge E$, using the
following homological argument.

\hfill

\lemma \label{_inje_on_H^2_calo_Lemma_}
Let $M$ be a compact K\"ahler manifold, $Z\subset M$
a real analytic subvariety of codimension at least 6.
Then the natural restriction
map $H^2(\calo_{M}) \arrow H^2(\calo_{M\backslash Z})$
is injective.

\hfill

{\bf Proof}
There is a natural map $\phi$ from the $p$-th de Rham cohomology
of a complex manifold $M$ to its $p$-th holomorphic cohomology
$H^p(M, \calo_M)$: given a closed $p$-form $\eta$,
the $(0,p)$-part of $\eta$ is $\bar\6$-closed and
represents a class in $H^p(M, \calo_M)$.
Consider the commutative diagram  
\begin{equation}\label{_commu_H^2_CD_Equation_}
\begin{CD}
H^2_{DR}(M) @>j>> H^2_{DR}(M\backslash Z)\\
@VV\phi V @VV {\phi_Z} V\\
H^2(\calo_{M}) @>{j_0}>> H^2(\calo_{M\backslash Z})
\end{CD}
\end{equation}
(here $H^*_{DR}$ denotes
the de Rham cohomology group). By definition,
$\ker \phi$, $\ker \phi_Z$ is the $F^1H^2_{DR}$-part
of $H^2_{DR}(M)$, $H^2_{DR}(M\backslash Z)$, where
$F^i$ denotes the Hodge filtration on cohomology.
The Dolbeault spectral sequence gives the $E_2$-term corresponding
to $F^1H^2_{DR}(V)$ for any complex manifold $V$ as follows:
\[ 0 \arrow H^0(\Omega^2(V))\arrow E_2(F^1H^2_{DR}(V)) \arrow 
     H^1(\Omega^1(V))\arrow 0.
\]
Using the functoriality of Dolbeault spectral sequence,
we obtain the following diagram with exact rows
\begin{equation}\label{_restri_E_2_Equation_}
\begin{array}{ccccccccc}
0 &\!\!\!\!\rightarrow\!\!\!\! & H^0(\Omega^2(M)) & \!\!\!\!\rightarrow\!\!\!\! & E_2(F^1H^2_{DR}(M)) &\!\!\!\!\rightarrow\!\!\!\! &
     H^1(\Omega^1(M)) & \!\!\!\!\rightarrow\!\!\!\! & 0\\
& & \downarrow & & \downarrow && \downarrow & & \\
0 &\!\!\!\!\rightarrow\!\!\!\!& H^0(\Omega^2(M\backslash Z))&\!\!\!\!\rightarrow\!\!\!\!& E_2(F^1H^2_{DR}(M\backslash Z)) &\!\!\!\!\rightarrow\!\!\!\!& 
     H^1(\Omega^1(M\backslash Z)) &\!\!\!\!\rightarrow\!\!\!\!& 0
\end{array}
\end{equation}
Using \ref{_H^i_coh_codim_Proposition_}, we obtain that the
vertical arrows of \eqref{_restri_E_2_Equation_} are isomorphisms.
Therefore, the $E_2(F^1H^2_{DR}$)-terms for $M$ and
$M\backslash Z$ are isomorphic. Since the Dolbeault
spectral sequence for $M$ degenerates in $E_2$, 
and $E_2(F^1H^2_{DR}(M)) = E_2(F^1H^2_{DR}(M\backslash Z))$,
all differentials $d_i$, $i>2$ for the
Dolbeault spectral sequence of $M\backslash Z$
vanish on $F^1H^2_{DR}(M\backslash Z)$.
Also, for $E_2(F^qH^p_{DR}(M))$ ($p, q\leq 1$),
the natural restriction map
\[ E_2(F^qH^p_{DR}(M))\arrow E_2(F^qH^p_{DR}(M\backslash Z))
\]
is an isomorphism, as follows from
\ref{_H^i_coh_codim_Proposition_},
hence the differentials $d_i$, $i>2$
vanish on the terms $E_2(F^qH^p_{DR}(M\backslash Z))$ ($p,
q\leq 1$) as well. This implies that the Dolbeault
spectral sequence degenerates in
$E_2(F^1H^2_{DR}(M\backslash Z))$, and the bottom
row of \eqref{_restri_E_2_Equation_} gives an exact sequence
\[ 0 \arrow H^0(\Omega^2(M\backslash Z))\arrow F^1H^2_{DR}(M\backslash Z) \arrow 
     H^1(\Omega^1(M\backslash Z)) \arrow 0.
\]
Applying \eqref{_restri_E_2_Equation_} again, we obtain that the 
restriction map induces an 
isomorphism 
\[ 
F^1H^2_{DR}(M\backslash Z) \cong F^1H^2_{DR}(M).
\]
In terms of \eqref{_commu_H^2_CD_Equation_} this is interpreted
as an isomorphism $\ker \phi=\ker \phi_Z$.
The left arrow of \eqref{_commu_H^2_CD_Equation_} is surjective
because $M$ is K\"ahler. An elementary diagram chasing
using surjectivity of $\phi$ and 
$\ker \phi=\ker \phi_Z$ implies that
$j_0$ is indeed injective. We proved 
\ref{_inje_on_H^2_calo_Lemma_}. \endproof

\hfill 

Return now to the proof of \ref{_degree_of_posi+inva_Theorem_}.
Consider $[\eta]$ as an element of $H^2_{DR}(M\backslash Z)$.
Then $\phi_Z([\eta])= [\eta^{0,2}_K]_{\bar\6}$, where
$\eta^{0,2}_K$ denotes the $\eta^{0,2}_K$-component
of $\eta$, and $[\cdot ]_{\bar\6}$ its Dolbeault class
in $H^2(\calo_{M\backslash Z})$. Now, $[\eta^{0,2}_K]_{\bar\6}$
belongs to the image of $j_0(H^2(\calo_{M}))$, because
$\eta^{0,2}_K= \overline{\eta^{2,0}_K}$ is $L^1$-integrable,
and the cohomology of currents are equal to cohomology of forms. 
\ref{_inje_on_H^2_calo_Lemma_} implies now that 
\[ \phi([\eta])=\phi_Z([\eta])=[\eta^{0,2}_K]_{\bar\6}.\]
Using \eqref{_degree_via_E_Equaion_}
we may compute $\deg_I[\eta]$ as an integral
\[ \deg_I[\eta]= \lambda_n^{-1}\int_M \eta_K^{2,0}\wedge E 
\]
(this makes sense, because $\eta_K^{2,0}$ is $L^1$-integrable).
Then
\begin{equation}\label{_deg_via_eta_+_Equation_} 
  \deg_I [\eta] = 2^{-(n-1)}\int_M \eta_K^{2,0}\wedge E =
  \int_M \eta_+ \wedge \omega_I^{2n-1}
\end{equation}
(the first equation holds by 
\eqref{_E_via_omega_Equation_}, and  the second one is implied by 
$(\eta_K^{2,0})^{1,1}_I= \eta_+$, which is essentially
a statement of \ref{_decompo_eta_w-K_SU(2)_Lemma_}).

Since $\eta_+$ is a positive 2-form, 
the integral \eqref{_deg_via_eta_+_Equation_}
is non-negative, and positive unless
$\eta_+=0$. We have reduced \ref{_degree_of_posi+inva_Theorem_} 
to \ref{_closed_K-positive_is_L^1_Proposition_}.
\endproof


\section{Positive $(2,0)$-forms on hyperk\"ahler manifolds}
\label{_integra_K-posit_Section_}


The purpose of this Section is to prove
\ref{_closed_K-positive_is_L^1_Proposition_}.
We deduce this result from the quaternionic version
of the classical Sibony's Lemma and Skoda-El Mir Theorem
(\cite{_El_Mir_}, \cite{_Skoda_}, \cite{_Sibony_}, \cite{_Demailly:L^2_}), 
which is proven in \cite{_Verbitsky:Skoda_}.

\hfill

\theorem\label{_Skoda_El_Mir_+_Sibony_Theorem_}
Let $(M,I,J,K,g)$ be a hyperk\"ahler manifold, and
$Z\subset (M,I)$ a closed complex subvariety. Consider a (2,0)-form
 $\eta$ on $(M,I)\backslash Z$, which satisfies the following conditions
\begin{description}
\item[(a)] $J\eta = \bar\eta$
\item[(b)] $\eta (x, J \bar x) \geq 0$, for any $x\in T^{1,0}(M,I)$.
\item[(c)] $\6_I\eta=0$.
\end{description}
Then
\begin{description}
\item[(i)] (Sibony's Lemma) The form $\eta$ is locally
integrable around $Z$, if $Z$ is compact and satisfies 
 $\codim_\C Z\geq 3$.
\item[(ii)] (Skoda-El Mir theorem)
Assume that $\eta$ is locally integrable around $Z$.
Consider the trivial extension $\tilde \eta$ of $\eta$ to
$M$ as a 
$(2, 0)$-current
on $(M,I)$. Then $\6_I\tilde\eta=0$.
\end{description}

{\bf Proof:} See \cite{_Verbitsky:Skoda_}, Theorem 1.2
and Theorem 1.3. \endproof

\hfill

Now we deduce 
\ref{_closed_K-positive_is_L^1_Proposition_} from
\ref{_Skoda_El_Mir_+_Sibony_Theorem_}.

\hfill

We work in assumptions of \ref{_closed_K-positive_is_L^1_Proposition_}.
Let $\rho \in\Lambda^{2,0}(M\backslash Z,K)$ be a $K$-positive, 
$\6_K$-closed form on $M\backslash Z$.
Since $SU(2)$ acts on the quaternionic triples $(I,J,K)$
transitively, there exists $g\in SU(2)$
mapping $K$ to $I$ and $I$ to $J$, with
$SU(2)$ acting on quaternions by conjugation. 
Let 
\[  \eta:= g(\rho)\in \Lambda^{2,0}(M\backslash Z,I)
\]
be the $(2,0)$-form on $(M\backslash Z,I)$ corresponding to 
$\rho$ under the isomorphism induced by $g$.

Since $I(\rho)= \bar\rho$, $\eta$ satisfies
$J(\eta)= \bar\eta$. Also, $K$-positivity
of $\rho$ is equivalent to
$\eta (x, I \bar x) \geq 0$, for any $x\in T^{1,0}(M,I)$.
To apply \ref{_Skoda_El_Mir_+_Sibony_Theorem_}
to $\eta$, it remains to show that $\6_I\eta=0$.

\hfill

Consider the complex vector space
\[ {\cal H}= \langle d, IdI, JdJ, KdK\rangle \subset \End(\Lambda^*(M)),\]
generated by the de Rham differential and its twists. Clearly,
${\cal H}$ is preserved by the natural action of $SU(2)$
on $\End(\Lambda^*(M))$. As  a representation of $SU(2)$,
${\cal H}$ has weight 1, that is, the standard
generator $\1 I$ of the Cartan subalgebra
of $SU(2)$ acts on ${\cal H}$ with eigenvalues
$\1$ and $-\1$. Therefore, ${\cal H}$
is isomorphic to a sum of two irreducible,
2-dimensional representations of $SU(2)$.
This implies that the space 
\[
{\cal H}^{1,0}_K:= \{ \delta \in {\cal H} \ \ | \ \
K(\delta) = \1\delta\}
\]
is 2-dimensional. This space is clearly
generated by $\6_K$ and $I \bar\6_K I$.
For any $\delta \in {\cal H}^{1,0}_K$,
$\delta(\rho)=0$, because $\6_K(\rho) = I \bar\6_K I(\rho)=0$.
Since 
\[ 
  g({\cal H}^{1,0}_K) = {\cal H}^{1,0}_I:= \{ \delta \in {\cal H} \ \ | \ \
I(\delta) = \1\delta\},
\]
all elements of ${\cal H}^{1,0}_I$ vanish on
$g(\rho)=\eta$. Therefore, $\6_I\eta=0$.
We obtain that $\eta$ satisfies conditions
(a), (b), (c) of \ref{_Skoda_El_Mir_+_Sibony_Theorem_}.

By \ref{_Skoda_El_Mir_+_Sibony_Theorem_} (i),
 $\eta$ is locally integrable; this proves 
\ref{_closed_K-positive_is_L^1_Proposition_} (i).
By \ref{_Skoda_El_Mir_+_Sibony_Theorem_} (ii),
a trivial extension of $\eta$ to $M$ is
$\6_I$-closed as a current. Backtracking the above
argument, we find that this is equivalent
to $\6_K$-closedness of the trivial extension
of $\rho$. We proved
 \ref{_closed_K-positive_is_L^1_Proposition_} (ii).


\section{Stability of weakly hyperholomorphic sheaves}
\label{_stabi_final_Section_}


In this section, we use \ref{_wea_hyhe_c_1_SU(2)_inv_Theorem_} 
and \ref{_degree_of_posi+inva_Theorem_} to prove
\ref{_wea_hh_hh_Theorem_}.

Let $M$ be a complex hyperk\"ahler manifold,
$I$ an induced complex structure, and  $F$
a reflexive sheaf on $(M,I)$. Assume that $F$
is weakly hyperholomorphic.

We are going to show that
$F$ is polystable. On contrary,
let $F'\subset F$ be a destabilizing subsheaf. Replacing
$F'$ by its reflexization if necessary, we may assume that
$F'$ is reflexive. 

Let $Z$ be the union of singular sets of $F$ and $F'$.
By \ref{_singu_codim_3_refle_Lemma_}, $\codim_\C Z\geq 3$.
Denote by $U$ the complement $U:= M\backslash Z$.
On $U$, both sheaves $F$ and $F'$ are bundles.
Let $\Theta\in \Lambda^{1,1}(U)\otimes \End(F)$ be the 
curvature of $F$, and  
$\Theta'\in \Lambda^{1,1}(U)\otimes \End(F')$ the 
curvature of $F'$. Denote by $A\in \Lambda^{1,0}(U, \Hom(F/F', F'))$
the so-called {\em second fundamental form of a sub-bundle $F'$}
(\cite{_Griffi_Harri_}). The curvature of $F'$ can be expressed
through $\Theta$ and $A$ as 
\begin{equation}\label{_seco_form_curv_Equation_}
\Theta' = \Theta\restrict{F'} - A\wedge A^\bot,
\end{equation}
where $A^\bot\in \Lambda^{1,0}(U, \Hom(F',F/F'))$
is the Hermitian adjoint of $A$, and the form
\[ \Theta\restrict{F'}\in \Lambda^{1,1}(U)\otimes \End(F')
\]
is obtained from $\Theta$ by the orthogonal projection
$\End(F)\arrow \End(F')$.

The following claim is quite elementary.

\hfill

\claim\label{_A_vani_splits_Claim_}
Let $F$ be a reflexive coherent sheaf over a complex manifold $X$,
and $F'\subset F$ a reflexive subsheaf. Assume that
$F$ is equipped with a Hermitian structure outside of singularities.
Assume, moreover, that the second fundamental from of $F'\subset F$ vanishes.
Then 
\begin{equation}\label{_decompo_F_refle_Equation_}
F\cong F' \oplus (F/F')^{**},
\end{equation}
where $(F/F')^{**} = \Hom(\Hom(F/F', \calo_X), \calo_X)$
denotes the reflexive hull of $F/F'$.

\hfill

{\bf Proof:} Let $Z$ be the union of singular sets of $F$ and $F'$.
By \ref{_singu_codim_3_refle_Lemma_}, $\codim_\C Z\geq 3$.
Consider the orthogonal decomposiion
\begin{equation}\label{_decompo_F_Equation_}
   F\restrict{X\backslash Z}\cong F'\restrict{X\backslash Z} 
   \oplus (F/F')\restrict{X\backslash Z},
\end{equation}
By the definition
of the second fundamental form, the connection on $F$
preserves the decomposition \eqref{_decompo_F_Equation_}
if and only if this form vanishes (see \cite{_Griffi_Harri_}).
Denote by $j:\; (X\backslash Z) \hookrightarrow X$ the natural
embedding. By \cite{_OSS_}, Ch. II, Lemma 1.1.12
(see \ref{_refle_obtained_Remark_}), we have
$j_*j^* F = F$, $j_*j^* F' = F'$, 
$j_*j^* (F/F') = (F/F')^{**}$.
Comparing this with \eqref{_decompo_F_Equation_},
we obtain the decomposition
\eqref{_decompo_F_refle_Equation_}.
\ref{_A_vani_splits_Claim_} is proven.
\endproof

\hfill

Return to the proof of \ref{_wea_hh_hh_Theorem_}.
By \ref{_wea_hyhe_c_1_SU(2)_inv_Theorem_},
the cohomology class $c_1(F)$ is $SU(2)$-invariant.
Therefore, by \ref{_Lambda_of_inva_forms_zero_Lemma_},
$\deg c_1(F)=0$. To show that $F$ is polystable, we need to show
that for any reflexive subsheaf $F'\subset F$, we have
$\deg c_1(F')\leq 0$, and if the equality is reached,
then the decomposition \eqref{_decompo_F_refle_Equation_}
holds.

By \eqref{_seco_form_curv_Equation_}, we have
\begin{equation}\label{_Tr_via_seco_form_Equation_}
-\Tr \Theta' = -\Tr\left(\Theta\restrict{F'}\right) + \Tr(A\wedge A^\bot)
\end{equation}
This form represents $-c_1(F')$. The 
first summand of the right hand side of 
\eqref{_Tr_via_seco_form_Equation_} is $SU(2)$-invariant.
Indeed, the form $\Theta$ is by our assumptions
$SU(2)$-invariant, and $\Theta\restrict{F'}$
is obtained by orthogonal projection from 
$\Lambda^2(M)\otimes \End (F)$ to  $\Lambda^2(M)\otimes \End (F')$,
but this projection obviously commutes with $SU(2)$.
The second summand of the right hand side of 
\eqref{_Tr_via_seco_form_Equation_}is manifestly positive. 
We arrive in the situation which is
close to that dealt in \ref{_degree_of_posi+inva_Theorem_}.
To apply \ref{_degree_of_posi+inva_Theorem_},
we need a closed 2-form $\eta$ on $U$ which is a sum of
an $SU(2)$-invariant form and a positive
form of weight 2; however, there is no
reason why $\Tr(A\wedge A^\bot)$ should be pure
of weight 2. Therefore, to use 
\ref{_degree_of_posi+inva_Theorem_},
we need the following elementary 
linear-algebraic lemma.

\hfill

\lemma\label{_decompo_of_posi_Lemma_}
Let $U$ be a hyperk\"ahler manifold,
$I$ an induced complex structure on $U$ and
$\eta$  a positive $(1,1)$-form on $(U,I)$.
Consider the decomposition
\[ \eta = \eta_0+\eta_+,\]
where $\eta_0$ is $SU(2)$-invariant, and
$\eta_+$ pure of weight 2. Then $\eta_+$ 
is positive; moreover, $\eta_+=0 \Rightarrow \eta=0$.

\hfill

{\bf Proof:} 
Let $K$ be an induced complex structure satisfying 
$I\circ K = - K \circ I$. Consider the multiplicative
action of $K$ on $\Lambda^*(M)$ defined in
\ref{_decompo_eta_w-K_SU(2)_Lemma_}. By
\ref{_decompo_eta_w-K_SU(2)_Lemma_},
$\eta_+ = \frac{1}{2}(\eta - K(\eta))$.
On the other hand, the cone of positive forms is generated
by the 2-forms 
\[ \eta_z:=\1 z\wedge \bar z,
\] 
where 
$z\in \Lambda^{1,0}(M,I)$. 
Clearly, \[ - K(\eta_z) =\1 K(z)\wedge  K(\bar z).\]
Since $K(\bar z)\in \Lambda^{1,0}(M,I)$, the form
$-K(\eta_z)$ is positive.

We have shown that, $\eta_+ = \frac{1}{2}(\eta - K(\eta))$ is a sum of 
two positive forms; this form is positive, and it is non-zero unless
$\eta=0$. This proves \ref{_decompo_of_posi_Lemma_}.\endproof

\hfill

Return to the situation described by
\eqref{_Tr_via_seco_form_Equation_}. We have shown that
the 2-form $(-\Tr \Theta')$ is a sum of a positive form 
$\mu= \Tr\left(A\wedge A^\bot\right)$ and an
$SU(2)$-invariant form. Decomposing $\mu$
onto $SU(2)$-invariant and pure weight 2 parts
as in \ref{_decompo_of_posi_Lemma_}, we find that
$-\Tr \Theta'$ is a sum of an $SU(2)$-invariant form
and a positive form $\eta_+$ which is pure of weight 2
with respect to the $SU(2)$-action.

Now we can apply \ref{_degree_of_posi+inva_Theorem_}.
We find that $\deg \Tr \Theta'\leq 0$, and the equality is reached
only if $\eta_+=0$. 
If $F'$ is destabilizing, we have $\deg \Tr \Theta'\geq 0$.
By  \ref{_decompo_of_posi_Lemma_}, this 
is equivalent to $\mu=0$, or, what is the same, $A=0$.
Now, if $A=0$, then $F$ splits as in 
\ref{_A_vani_splits_Claim_}.

We have shown that for any destabilizing 
reflexive subsheaf $F'\subset F$,
the sheaf $F$ splits as $F= F' \oplus (F/F')^{**}$.
This means that $F$ is polystable.

\hfill

\remark\label{_SU(2)-inv_c_i_Remark_}
In assumptions of \ref{_wea_hh_hh_Theorem_},
one can easily see that all stable direct
summands $F_i$ of $F$ are also hyperholomorphic,
that is, have $SU(2)$-invariant $c_1$, $c_2$.
Indeed, $c_1(F_i)$ can be computed using the
curvature of $\nabla$ as indicated above, 
and by \eqref{_Tr_via_seco_form_Equation_}
it is $SU(2)$-invariant. For any coherent sheaf $A$,
let \[ D(A):= 2c_2(A) - \frac{\rk A-1}{\rk A}c_1(A)^2\]
be the {\em discriminant} of $A$. 
For any stable coherent sheaf $A$ with $SU(2)$-invariant
first Chern class, $D(A)$ satisfies the inequality 
\begin{equation}\label{_discri_inequa_Equation_}
\int_M \omega_J^{2n-2}\wedge D(A) \leq \int_M \omega_I^{2n-2}\wedge  D(A)
\end{equation}
(\cite{_V:Hyperholo_sheaves_}, Claim 3.21),
and the equality is reached of and only $A$ is
hyperholomorphic.
Since the discriminant is additive, and
\eqref{_discri_inequa_Equation_} is true for the
direct sum of all $F_i$, it is true for each summand
$F_i$. 

\hfill

\hfill

{\bf Acknowledgements:}
Many thanks  to S. Alesker, D. Kaledin  and the referee 
for useful suggestions and finding gaps in my arguments,
S. Merkulov for his interest and encouragement, 
D. Alekseevsky, R. Bielawski, F. Bogomolov, and T. Pantev for 
interesting discussions, and J.-P. Demailly
for answering my questions about positive forms. 
And last but not least, my gratitude
to Glasgow University, where I visited 
while writing this paper, for warmth 
and hospitality.

\hfill

{\small

\hfill

\noindent {\sc Misha Verbitsky\\
Laboratory of Algebraic Geometry, SU-HSE, \\
7 Vavilova Str., Moscow, Russia, 117312
 }\\
\tt  verbit@mccme.ru 
}


\begin{thebibliography}{6666}

\bibitem[AV]{_Alesker_Verbitsky_HKT_} 
Semyon Alesker, Misha Verbitsky,
{\em Plurisubharmonic functions on hypercomplex manifolds and
 HKT-geometry,} J. Geom. Anal. 16 (2006), no. 3, 375--399.


\bibitem[BS]{_Bando_Siu_} 
 Bando, S., Siu, Y.-T, 
{\it Stable sheaves and Einstein-Hermitian metrics}, 
In: Geometry and Analysis on
Complex Manifolds, Festschrift for Professor S. Kobayashi's 60th Birthday,
ed. T. Mabuchi, J. Noguchi, T. Ochiai, World Scientific, 1994, pp. 39-50.


\bibitem[Bea]{_Beauville_} 
 Beauville, A. {\em 
Varietes K\"ahleriennes dont la premi\`ere classe de Chern est
nulle.}  J. Diff. Geom. {\bf 18}, pp. 755-782 (1983).


\bibitem[Bes]{_Besse:Einst_Manifo_} 
Besse, 
A., {\em Einstein Manifolds}, Springer-Verlag, New York (1987)

\bibitem[BBR]{_BBR_3} 
Bartocci, C., Bruzzo, U., Hernandez Ruiperez, D.
{\em A hyper-K\"ahler Fourier transform}
Differential Geom. Appl. 8 (1998), no. 3, 239--249. 

\bibitem[Bo]{_Bogomolov_} Bogomolov, F.A., {\em Hamiltonian K\"ahler
manifolds}, Sov. Math. Dokl. {\bf 19} (1978), 1462--1465.

\bibitem[Br]{_Bryant_holo_87_}
Bryant, R., {\em  Metrics with exceptional holonomy}, 
Ann. of Math. 126 (1987), 525-576.



\bibitem[Ca]{_Calabi_}  Calabi,  E.,
{\em Metriques k\"ahleriennes et fibr\`es holomorphes}, 
Ann. Ecol. Norm. Sup. {\bf 12} (1979), 269-294.  



\bibitem[D1]{_Demailly_} 
Demailly, J.-P.,
{\em Complex analytic and algebraic geometry},
a book, \\{\small http://www-fourier.ujf-grenoble.fr/~demailly/books.html}

\bibitem[D2]{_Demailly:L^2_}
Demailly, Jean-Pierre, {\em  $L^2$ vanishing theorems for
  positive line bundles and adjunction theory},
Lecture Notes of a CIME course on "Transcendental Methods
of Algebraic Geometry" (Cetraro, Italy, July 1994),  arXiv:alg-geom/9410022,
and also Lecture Notes in Math., 1646, pp. 1--97, Springer, Berlin, 1996

\bibitem[E]{_El_Mir_}
H. El Mir, {\em Sur le prolongement des courants positifs
fermes,} Acta Math., 153 (1984), 1-45.



\bibitem[GH]{_Griffi_Harri_} 
Griffiths, Ph., Harris, J., 
{\em Principles of Algebraic Geometry}, \\ Wiley-Intersience,
New York, 1978.


\bibitem[KV]{_NHYM_} 
Kaledin, D., Verbitsky, M.,
{\it Non-Hermitian 
Yang-Mills connections},
alg-geom 9606019 (1996), 48 pages, LaTeX 2e.
(also published in Selecta Math. (N.S.) {\bf 4}
(1998), no. 2, 279--320)



\bibitem[Kob]{_Kobayashi_} 
Kobayashi S., {\em Differential geometry of complex vector bundles,} 
Princeton University Press, 1987.

\bibitem[OSS]{_OSS_} 
 Christian Okonek, Michael 
Schneider, Heinz  Spindler,
{\it Vector bundles on complex projective spaces.}
 Progress in mathematics, vol. 3,
 Birkhauser, 1980.

\bibitem[Sa]{_Salamon_quate_} 
Salamon, Simon M.
{\em Differential geometry of quaternionic manifolds}, 
Annales Scientifiques de l'Ecole Normale Suprieure 
Ser. 4, 19 no. 1 (1986), p. 31-55.


\bibitem[Sch]{_Scheja_} 
Scheja, G., {\it Riemannische Hebbarkeitss\"atze
f\"ur Cohomologieklassen}, Math. Ann. 144 (1961) pp. 345-360

\bibitem[Sib]{_Sibony_}
Sibony, Nessim,
{\em Quelques problemes de prolongement de courants en analyse complexe,}
Duke Math. J. 52, 157-197 (1985).

\bibitem[Sim]{_Simpson_} 
C.T. Simpson, 
{\em Constructing variations of Hodge structure using Yang-Mills
theory, and applications to uniformization}, 
Jour. of Amer. Math. Soc. {\bf 4} (1988), 867-918. 

\bibitem[Sk]{_Skoda_}
H. Skoda,
{\em Prolongement des courants positifs fermes de masse finie,} 
Invent. Math., 66 (1982), 361-376.

\bibitem[V0]{_Verbitsky:Symplectic_II_} 
Verbitsky M., {\em Hyperk\"ahler embeddings and holomorphic 
symplectic geometry II,} alg-geom electronic preprint 9403006 (1994),
14 pages, LaTeX,
also published in: GAFA {\bf 5} no. 1 (1995), 92-104.


\bibitem[V1]{_Verbitsky:Hyperholo_bundles_} 
Verbitsky M., 
{\em Hyperholomorphic bundles over a hyperk\"ahler manifold}, 
alg-geom electronic preprint 9307008 (1993), 43 pages, LaTeX,\\
also published in: 
Journ. of Alg. Geom., {\bf 5} no. 4 (1996) pp. 633-669.

\bibitem[V2]{_Verbitsky:Symplectic_I_} 
{\em Hyperk\"ahler and holomorphic 
         symplectic geometry I,} 
arXiv:alg-geom/9307009, 
Journ. of Alg. Geom., vol. 5 no. 3 pp. 401-415 (1996).


\bibitem[V3]{_V:Hyperholo_sheaves_} 
Verbitsky M., {\em
Hyperholomorphic sheaves and new examples of hyperk\"ahler manifolds,}
alg-geom 9712012, published in a book 
{\em ``Hyperkahler manifolds''}, D. Kaledin, M. Verbitsky,
International Press, Somerville, MA, 1999.

\bibitem[V4]{_Verbitsky:hypercomple_}
Verbitsky M.,
{\em Hypercomplex Varieties}, 
alg-geom/9703016 (1997); published in:
 Comm. Anal. Geom. {\bf 7} 
(1999), no. 2, 355--396.

\bibitem[V5]{_Verbitsky_HKT-exa_}
M. Verbitsky,
{\em Hyperk\"ahler manifolds with torsion obtained from hyperholomorphic bundles}
 math.DG/0303129, (Math. Res. Lett. 10 (2003), no. 4, 501--513).

\bibitem[V6]{_Verbitsky:qD_}
M. Verbitsky,
{\em Quaternionic Dolbeault complex and vanishing theorems on
hyperkahler manifolds}, math/0604303, 
Compos. Math. 143 (2007), no. 6, 1576--1592.


\bibitem[V7]{_Verbitsky:omega-psh_}
M. Verbitsky,
{\em Plurisubharmonic functions in calibrated geometry and q-convexity},
 arXiv:0712.4036, Math. Z., Vol. 264, No. 4, pp. 939-957 (2010).

\bibitem[V8]{_Verbitsky:Skoda_}
M. Verbitsky,
{\em Positive forms on hyperkahler manifolds,}
 arXiv:0801.1899, Osaka J. Math. Volume 47, Number 2 (2010), 353-384.

\bibitem[UY]{_Uhle_Yau_}
 Uhlenbeck K., Yau S. T.,  {\em On the existence of
Hermitian Yang-Mills connections in  stable vector bundles}, 
Comm. on Pure and Appl. Math., 
{\bf 39}, p. S257-S293 (1986).


\bibitem[Y]{_Yau:Calabi-Yau_} 
Yau, S. T., {\em On the Ricci curvature of a compact K\"ahler manifold 
and the complex Monge-Amp\`ere equation I.} Comm. on Pure and Appl.
Math. 31, 339-411 (1978).




\end{thebibliography}
\end{document}